\documentclass{article}
\usepackage{epsfig,theorem,amssymb}

\newtheorem{lemma}{Lemma}
\newtheorem{conj}{Conjecture}
\newtheorem{cor}{Corollary}
\newtheorem{ta}{Table}
\newtheorem{theorem}{Theorem}

\theoremstyle{plain}
\theorembodyfont\rm
\newtheorem{example}{Example}
\newtheorem{remark}{Remark}

\newcommand{\nc}{\newcommand}
\nc{\be}{\begin{equation}}
\nc{\ee}{\end{equation}}
\nc{\bea}{\begin{eqnarray}}
\nc{\eea}{\end{eqnarray}}
\nc{\disp}{\displaystyle}
\nc{\calN}{{\cal N}}
\nc{\calC}{{\cal C}}
\nc{\calM}{{\cal M}}
\nc{\calS}{{\cal S}}
\nc{\phit}{\hat{\varphi}}
\nc{\chit}{\hat{\chi}}
\nc{\hcalN}{\hat{\calN}}
\nc{\hcalS}{\hat{\calS}}
\nc{\hS}{\hat{S}}
\nc{\sigmad}{\sigma^\dagger}
\nc{\psid}{\psi^\dagger}

\font\tenmsb=msbm10
\font\sevenmsb=msbm7
\font\fivemsb=msbm5
\newfam\msbfam
\textfont\msbfam=\tenmsb
\scriptfont\msbfam=\sevenmsb
\scriptscriptfont\msbfam=\fivemsb

\def\sz{\scriptsize}
\def\sstyle{\scriptstyle}

\def\i{{\rm i}}
\def\F{{\mathcal F}}
\def\Sp{{\rm Span}}
\def\binom#1#2{{#1\choose #2}}

\begin{document}

\title{Loops, matchings and alternating-sign matrices
\footnote{Extended version of a talk given at the 14th International
Conference on Formal Power Series and Algebraic Combinatorics held
from 8-12 July 2002 in Melbourne, Australia.}}

\author{Jan de Gier\\
Department of Mathematics and Statistics,\\
The University of Melbourne,
Parkville, Victoria 3010, Australia}

\date{\today}

\maketitle

\begin{abstract}
The appearance of numbers enumerating alternating sign matrices in
stationary states of certain stochastic processes on matchings is
reviewed. New conjectures concerning nest distribution functions are
presented as well as a bijection between certain classes of
alternating sign matrices and lozenge tilings of hexagons with cut off
corners.
\end{abstract}

\section{Introduction}

Following an observation of Razumov and Stroganov \cite{RazuS01a} for
the XXZ spin chain (most terms will be defined below), it was observed
by Batchelor et al. \cite{BatchGN01} that the groundstates of loop
Hamiltonians with different boundary conditions are related to
symmetry classes of the alternating sign matrices introduced by Mills
et al. \cite{MillsRR82,MillsRR83,Robbins00}. The following is an account of
this surprising new relation between physics and combinatorics. 

The connection can be described using the action of the
Temperley-Lieb algebra on matchings of $\{1,\ldots,n\}$. This action,
which will be described in detail in Section \ref{se:Stoch}, will be
used to define the dense O(1) loop model. In Section \ref{se:FPL} it is
shown that matchings define an equivalence relation on alternating
sign matrices. These two ingredients have led several authors
\cite{MitraNGB02,PearceRGN02,RazuS01b,RazuS01c} to formulate
conjectures relating symmetry classes of alternating sign
matrices to boundary conditions in the dense O(1) loop model. The
exact correspondence is stated at the end of Section \ref{se:FPL}.

Many of the alternating sign matrix numbers factorise into small
primes. The appearance of such numbers in stationary states has a nice
application. As shown in Section \ref{se:Nests}, exact closed form
formulae of certain expectation values and correlation functions can
be guessed from exact calculations on a few examples. This is
extremely useful since one obtains exact conjectural results
for physical quantities. In fact, in some cases one may argue the
validity of the conjectured formulae, or at least their asymptotics,
because they make sense physically \cite{GierNPR03}. Identification of numbers can also
be useful mathematically, since relations between different
combinatorial objects may be discovered. For example, a direct
bijection between a class of alternating sign matrices and hexagons
with cut off corners was discovered this way. This bijection is
described in Section \ref{se:Hex}.

\section{A stochastic process on matchings}
\label{se:Stoch}
In the first two subsections we define a stochastic process on
non-crossing perfect matchings \cite{PearceRGN02}. A simple example of
such a process is given at the end of subsection \ref{se:ham}. In the
remaining subsections similar stochastic processes on more general
matchings will be defined.

\subsection{Matchings and matchmakers}
\label{se:match}
A {\bf $p$-matching}, or simply matching, of the {\bf vertex set}
$[n]=\{1,\ldots,n\}$ is an unordered collection of $p$ pairs of
vertices, or {\bf edges}, and $n-2p$ single vertices. A matching $F$
is called {\bf crossing} if it contains an edge $\{i,j\}$ and a vertex
$k$ such that $i<k<j$ or if it contains edges $\{i,j\}$ and $\{k,l\}$
such that $i<k<j<l$. A matching is {\bf perfect} if $n=2p$ and {\bf 
near-perfect} if $n=2p+1$. Let $\F_{2n}$ denote the set of
all non-crossing perfect matchings of $[2n]$, and $\F_{2n+1}$ the set
of all non-crossing near-perfect matchings of $[2n+1]$. 

\begin{example}
\label{ex:basis}
For $n=6$ there are five non-crossing perfect matchings.
\begin{eqnarray*}
\F_6 &=&
\{\{1,2\}\{3,4\}\{5,6\},\{1,4\}\{2,3\}\{5,6\},\{1,2\}\{3,6\}\{4,5\},\\
&&\hphantom{\{} \{1,6\}\{2,3\}\{4,5\},\{1,6\}\{2,5\}\{3,4\}\}.
\end{eqnarray*}
The edges of each matching are written here in a particular
order but it is to be understood that no particular order is
preferred. A notation for the matchings that will be useful
later is to depict them as loop segments connecting vertices. The five
matchings comprising $\F_6$ will thus be denoted by the following
pictures,
\[
\begin{picture}(240,80)
\put(0,0){\epsfxsize=240pt\epsfbox{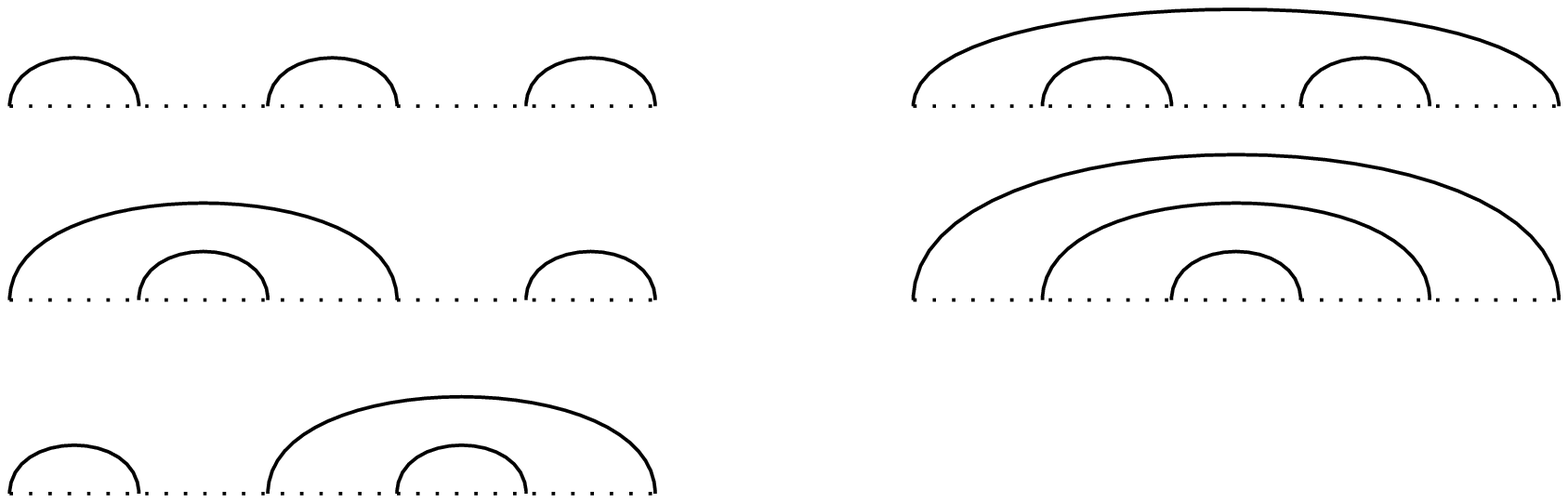}}
\put(-20,60){$1:$}
\put(-20,30){$2:$}
\put(-20,0){$3:$}
\put(120,60){$4:$}
\put(120,30){$5:$}
\end{picture}
\]
Instead of this graphical notation, a sometimes more convenient typographical
notation for matchings of $[n]$ is obtained by using parentheses for
paired vertices \cite{MitraNGB02},  
\[
\F_{6}\; =\; \{()()(),(())(),()(()),(()()),((()))\}.
\]
We will use the graphical and parenthesis notation interchangebly.
\end{example}

Define matching generators or {\bf matchmakers} $e_j$,
$j\in\{1,\ldots,2n-1\}$ acting non-trivially on elements $F\in\F_{2n}$
containing $j$ and $j+1$, and as the identity otherwise. With the
identification $\{i,k\}=\{k,i\}$ the vertices $j$ and $j+1$ can occur
in edges of $F$ in essentially two distinct cases. The action of $e_j$
in those cases is defined by     
\begin{equation}
e_j: \left\{ 
\begin{array}{lcl}
\{j,j+1\} & \mapsto & \{j,j+1\} \\
\{i,j\}\{j+1,k\} & \mapsto & \{i,k\}\{j,j+1\}\\
\end{array}\right.
\label{eq:TLgen_def}
\end{equation}
Equation (\ref{eq:TLgen_def}) defines the action of $e_j$ for
all orderings of $i,j,k$ with $j\in \{1,\ldots,2n-1\}$ and $i,k \in
\{1,\ldots,2n\}$. 
\begin{lemma}
The matchmakers $e_j$, $j\in\{1,\dots,2n-1\}$ satisfy the following relations,
\begin{eqnarray}
&&e_j^2 = (q+q^{-1}) e_j \nonumber\\
&&e_je_{j\pm1}e_j = e_j \label{eq:TLdef}\\
&&e_je_k = e_ke_j \quad |j-k| > 1, \nonumber
\end{eqnarray}
with $q=\exp(\i\pi/3)$. 
\end{lemma}

For general $q$ the algebra (\ref{eq:TLdef}) is called the {\bf
Temperley-Lieb algebra} \cite{TempL71}. There exists a graphical
representation of the $e_j$ which is closely related to the
graphical notation of matchings of $[2n]$,
\[
e_j \quad =\quad
\begin{picture}(140,15)
\put(2,-10){\epsfxsize=135pt\epsfbox{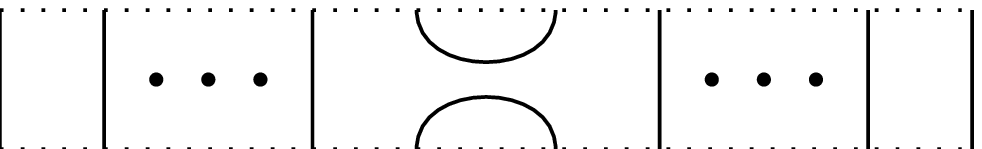}}
\put(.5,-18){\sz 1}
\put(15,-18){\sz 2}
\put(38,-18){$\sstyle j-1$}
\put(58,-18){$\sstyle j$}
\put(72,-18){$\sstyle j+1$}
\put(89,-18){$\sstyle j+2$}
\put(112,-18){$\sstyle 2n-1$}
\put(135.5,-18){$\sstyle 2n$}
\end{picture}\hspace{10pt}
\label{eq:monoid}
\]
\vskip20pt
\noindent
The graph of the multiplication $w_1w_2$ of two words in the
Temperley-Lieb algebra is obtained by placing the graph of $w_1$ below
the graph of $w_2$ and erasing the intermediate dashed line. The
algebraic relations (\ref{eq:TLdef}) now have a nice pictorial
interpretation.
\begin{example}
The relations $e_j^2=e_j$ and $e_je_{j+1}e_j=e_j$ are graphically depicted 
as
\[
\begin{picture}(110,40)
\put(0,0){\epsfxsize=40pt\epsfbox{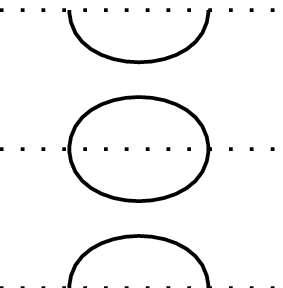}}
\put(50,18){$=$}
\put(70,10){\epsfxsize=40pt\epsfbox{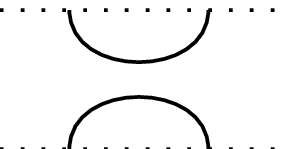}}
\end{picture}
\]
\vskip2mm
\[
\begin{picture}(150,60)
\put(0,0){\epsfxsize=60pt\epsfbox{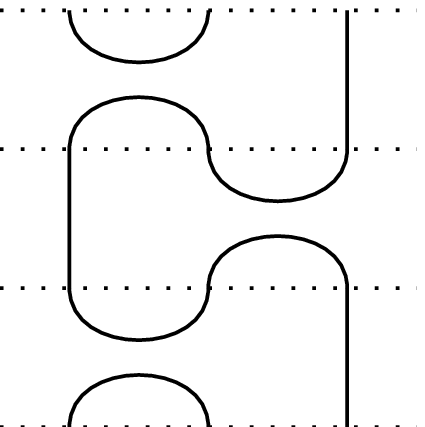}}
\put(70,28){$=$}
\put(90,20){\epsfxsize=60pt\epsfbox{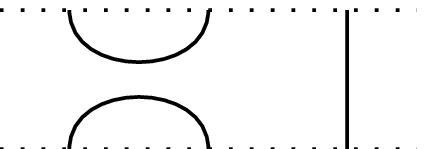}}
\end{picture}
\]
The action of $e_1$ on $(()())$ is given by\[
\begin{picture}(270,40)
\put(0,0){\epsfxsize=120pt\epsfbox{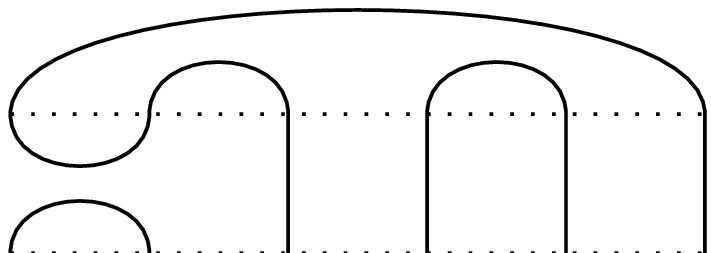}}
\put(130,10){$=$}
\put(150,10){\epsfxsize=120pt\epsfbox{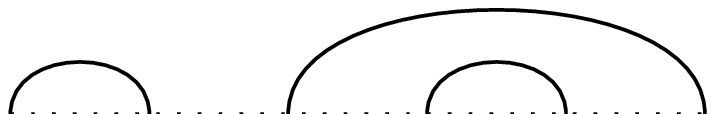}}
\end{picture}
\]
\end{example}

\subsection{Hamiltonian and stationary state}
\label{se:ham}
The loops in the graphical representation of the Temperley-Lieb
algebra have the physical interpretation of boundaries of percolation
clusters. A much studied object in physics is the loop energy
operator, or loop {\bf Hamiltonian}, which is defined as  
\begin{equation}
H^{\rm C}_{2n} = \sum_{j=1}^{2n-1} (1-e_j) \label{eq:hamCdef}
\end{equation}
The superscript C stands for {\bf closed boundary conditions}, which will be
explained in section \ref{se:bc}. In the representation $H^{\rm C}_{2n}:
\Sp(\F_{2n})\rightarrow \Sp(\F_{2n})$, $H^{\rm C}_{2n}$ is called the 
Hamiltonian of the {\bf Temperley-Lieb loop model}, or {\bf dense O(1)
loop model} with closed boundary conditions. We will refer to this
representation as the {\bf loop representation} of $H^{\rm C}_{2n}$.

In the loop representation $H^{\rm C}_{2n}$ is a matrix whose off-diagonal
entries are all non-positive and whose columns add up to zero
\cite{PearceRGN02}. Such a matrix is called an {\bf intensity matrix}
and defines a {\bf stochastic process} in continuous time given by the
master equation,
\[
\frac{{\rm d}}{{\rm d}t} P_{2n}(t) = -H^{\rm C}_{2n} P_{2n}(t),\qquad P_{2n}(t)
=\sum_{F\in\F_{2n}} a_F(t) F,
\label{eq:master}
\]
where $a_F(t)$ is the unnormalized probability to find
the system in the state $F$ at time $t$. Since $H^{\rm C}_{2n}$ is an intensity
matrix it has at least one zero eigenvalue. Its corresponding left
eigenvector $P_{2n}^{\rm L}$ is trivial and its right eigenvector
$P_{2n}$ is called the {\bf stationary state}, 
\begin{eqnarray*}
P_{2n}^{\rm L} H^{\rm C}_{2n} = 0,&&\qquad P_{2n}^{\rm L} = (1,1,\ldots,1),\\
H^{\rm C}_{2n} P_{2n} = 0,&&\qquad P_{2n} = \lim_{t\to\infty} P_{2n}(t).
\label{eq:stat}
\end{eqnarray*}
Properties of the stationary state from a physical perspective are
described in \cite{GierNPR02}.
\begin{example}
\label{ex:Cn=6}
For $n=3$ the action of $H^{\rm C}_6$ on $\Sp(\F_{6})$ can be
calculated by its action on the five basis states
()()(), (())(), ()(()), (()()), ((())) (see also Example \ref{ex:basis}).
We find for example,
\[
H^{\rm C}_6 ()()() = 2()()() - (())() - ()(()).
\]
Similarly calculating the action of $H^{\rm C}_6$ on the other basis states
yields
\[
H^{\rm C}_6 = - \left( \begin{array}{@{}ccccc@{}} 
-2 & 2 & 2 & 0 & 2 \\
1 & -3 & 0 & 1 & 0 \\
1 & 0 & -3 & 1 & 0 \\
0 & 1 & 1 & -3 & 2 \\
0 & 0 & 0 & 1 & -4 
\end{array}\right).
\]
The stationary state $P_6$ of $H^{\rm C}_6$ is given by
\[
P_6^{\rm T} = (11,5,5,4,1),
\]
where T denotes transposition.
\end{example}

The stationary state turns out to have a surprising
combinatorial interpretation. Let us look at a few more explicit
solutions of $H^{\rm C}_{2n}P_{2n}=0$,
\[
\begin{array}{c|c}
2n & P_{2n}^{\rm T} \\ \hline
2 & (1)\\
4 & (2,1)\\
6 & (11,5,5,4,1)\\
8 & (170,75,75,71,56,56,50,30,14,14,14,14,6,1)\\
\end{array}
\]
We can now make the surprising observation \cite{BatchGN01} that the
largest components of $P_{2n}$, i.e. $\{1,2,11,170,..\}$, enumerate
cyclically symmetric transpose complement plane partitions, and that
$P_{2n}^{\rm L} P_{2n}$ which is the sum of elements of $P_{2n}$,
i.e. $\{1,3,26,646...\}$, enumerate vertically symmetric alternating
sign matrices. 

{\bf Alternating sign matrices} (ASMs) were introduced by Mills et
al. \cite{MillsRR82,MillsRR83} and are matrices with 
entries in $\{-1,0,1\}$ such that the entries in each column and 
each row add up to $1$ and the non-zero entries alternate in sign.
Mills et al. conjectured that the number of ASMs is given by
\begin{equation}
A_n = \prod_{j=0}^{n-1} \frac{(3j+1)!}{(n+j)!} = 1,2,7,42,429,\ldots,
\label{eq:ASMno}
\end{equation}
which was proved more than a decade later by Zeilberger \cite{Zeilb96a} and
Kuperberg \cite{Kupe96}. Conjectured enumerations of symmetry classes were given 
by Robbins \cite{Robbins00}, many of which were subsequently proved by
Kuperberg \cite{Kupe00}. Proofs of some remaining conjectures were
announced recently by Okada \cite{Odada02}. The properties and
history of ASMs are  reviewed in a book by Bressoud \cite{Bress99}, as
well as by Robbins \cite{Robbins91} and Propp \cite{Propp01}.

In the following we will see that not only the largest components and
the sum of components of the stationary state have a meaning, but that
the other integers also have a combinatorial interpretation
\cite{RazuS01b}. It turns out that other symmetry classes of ASMs
appear when we use different boundary conditions for the loop
Hamiltonian. To show that, we first need to extend some definitions
and introduce a few more concepts. 

\begin{remark}
\label{rem:Dyck}
It is sometimes useful to think of the
non-crossing perfect matchings as Dyck paths. In terms of the
parentheses notation, a Dyck path is obtained from each non-crossing
matching by moving a step in the NE direction for each opening
parenthesis `$($' and a step SE for each closing parenthesis `$)$', or
in terms of pictures, 
\[
\begin{picture}(270,30)
\put(0,0){$()(())\; \sim$}
\put(45,0){\epsfxsize=120pt\epsfbox{e1_on_w2.eps}}
\put(170,0){$\sim$}
\put(190,0){\epsfxsize=80pt\epsfbox{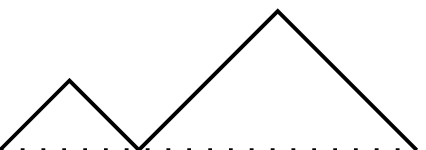}}.
\end{picture}
\]
\end{remark}

\begin{remark}
There exists a representation of the Temperley-Lieb algebra on the space
$V = \otimes_{i=1}^{2n} {\mathbb C}^2$ where $e_i$ is represented by the
following $4\times4$ matrix on the $i$th and $(i+1)$th copy of
${\mathbb C}^2$ and as the identity elsewhere,
\[
e = \left(\begin{array}{@{}cccc@{}}
0 & 0 & 0 & 0\\
0 & q & -1 & 0 \\
0 & -1 & q^{-1} & 0\\
0 & 0 & 0 & 0
\end{array}\right),
\]
with $q=\exp (\i\pi/3)$. In this representation $H^{\rm C}_{2n}$ is called the
Hamiltonian of the XXZ spin chain at $\Delta= -(q+q^{-1})/2 = -1/2$
with diagonal open boundary conditions, which is closely related to
the six-vertex model, see e.g. \cite{Baxter82,KoreIB93}.
\end{remark}

\begin{remark}
The Temperley-Lieb algebra is also closely related to the study of knot
invariants and the Jones polynomial, see e.g. \cite{Jones85,Kauff87}.
\end{remark}

\subsection{Boundary conditions}
\label{se:bc}
A matching is {\bf directed} if we distinguish between
$\{i,j\}$ and $\{j,i\}$. Directed edges will be denoted by
$(i,j)$. Unless stated otherwise, matchings will be non-directed. 
A {\bf left extended} $(p,k)$-matching of $[n]$ is
obtained from a $p$-matching by pairing $k$ unmatched vertices with
an additional vertex labelled $0$. A {\bf right extended}
$(p,k)$-matching of $[n]$ is obtained from a $p$-matching by pairing
$k$ unmatched vertices with and additional vertex labelled $n+1$. An
{\bf extended} $(p,k_1,k_2)$-matching of $[n]$ is a left extended
$(p,k_1)$-matching and a right extended $(p,k_2)$-matching. An
extended $(p,k_1,k_2)$-matching of $[n]$ is perfect if
$2p+k_1+k_2=n$. Let $\F^{\rm re}_n$ the set of all non-crossing
perfect right extended matchings of $[n]$ and $\F^{\rm e}_n$ the set
of all non-crossing perfect extended matchings of $[n]$. The
parentheses notation carries over to extended matchings in an obvious
way.

\begin{remark}
A $p$-matching can be identified with a perfect left or right extended
$(p,n-2p)$-matching.
\end{remark}

\begin{example}
For $n=4$ there are six non-crossing perfect right extended matchings,
two of which are perfect matchings,
\[
\F^{\rm re}_{4} = \{((((,(((),(()(,()((,()(),(())\},
\]
and there are six non-crossing perfect extended matchings for $n=3$,
\[
\F^{\rm e}_3=\{))),()), ()(,((),)(),(((\}.
\] 
\end{example}

The action of the generators $e_j,\; j\in\{1,\ldots,n-1\}$ defined in
(\ref{eq:TLgen_def}) carries over to extended matchings where now $i,k\in
\{0,1,\ldots,n,n+1\}$, and with the additional definitions  
\[
e_j: \left\{ 
\begin{array}{lclc}
\{j,n+1\}\{j+1,n+1\} & \mapsto & \{j,j+1\} \\
\{0,j\}\{0,j+1\} &\mapsto & \{j,j+1\}\\
\{0,j\}\{j+1,n+1\} &\mapsto & \{j,j+1\}
\end{array}\right.
\]
In the pictorial representation these relations allow us to erase
those parts of a picture that are connected to the external vertices but
are disconnected from the vertex set $[n]$. The action of $e_i$ on
$F\in\F^{\rm e}_{n}$ can again be obtained graphically by placing the
graph of $e_i$ below that of $F$ and erasing all disconnected
parts.  
\begin{example}
\label{ex:bound}
$e_5 (())(( = (())() \in \F^{\rm e}_6$ can be displayed graphically as
\[
\begin{picture}(270,40)
\put(0,0){\epsfxsize=120pt\epsfbox{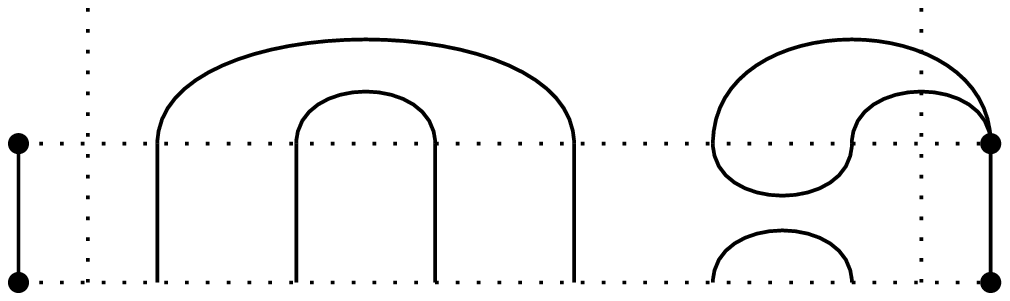}}
\put(130,10){$=$}
\put(150,10){\epsfxsize=120pt\epsfbox{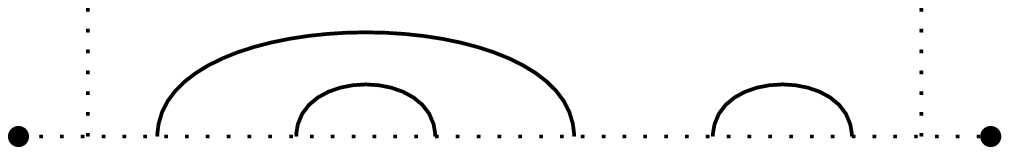}}
\end{picture}
\]
The external vertices are represented by bold dots. In the following we will
only draw the part inside the boundaries denoted by the vertical lines.
\end{example}

\subsubsection{Open and mixed boundary conditions}
In \cite{MitraNGB02} {\bf boundary generators} $f_1, f_n$ are introduced
that act non-trivially on elements of $F\in\F^{\rm e}_{n}$ containing
$1$ or $n$ respectively, 
\[
f_1: \left\{ 
\begin{array}{lcl}
\{0,1\} & \mapsto & \{0,1\} \\
\{1,n+1\} & \mapsto & \{0,1\} \\
\{1,i\} & \mapsto & \{0,1\}\{0,i\}
\end{array}\right.
\]
with $i\in\{2,\ldots,n\}$, and
\[
f_n: \left\{ 
\begin{array}{lcl}
\{n,n+1\} & \mapsto & \{n,n+1\} \\ 
\{0,n\} & \mapsto & \{n,n+1\}\\
\{i,n\} & \mapsto & \{i,n+1\}\{n,n+1\} \\
\end{array}\right.
\]
with $i\in\{1,\ldots,n-1\}$.

\begin{lemma}
The generators $f_1$ and $f_n$ satisfy the relations
\begin{eqnarray*}
&&f_1^2 = f_1,\qquad e_1 f_1 e_1 = e_1,\\
&&f_n^2 = f_n,\qquad e_{n-1} f_n e_{n-1} = e_{n-1},
\end{eqnarray*}
They are graphically represented by
\[
\begin{picture}(240,20)
\put(0,10){$f_1=$}
\put(30,0){\epsfxsize=70pt\epsfbox{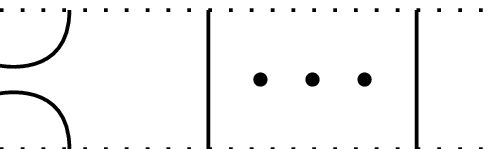}}
\put(140,10){$f_n=$}
\put(170,0){\epsfxsize=70pt\epsfbox{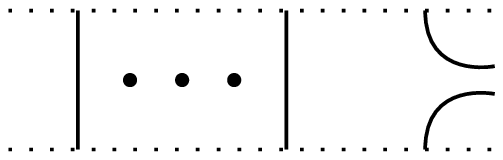}}
\end{picture}
\]
\end{lemma}
Furthermore, the definition of $f_1$ and $f_n$ is such that we can
erase those parts of composite pictures that connect the left and
right external vertices but are otherwise disconnected.
\begin{example}
\label{ex:Obc-rem}
$e_5 (()))( = (())() \in \F^{\rm e}_6$ can be displayed graphically as 
\[
\begin{picture}(270,50)
\put(0,0){\epsfxsize=120pt\epsfbox{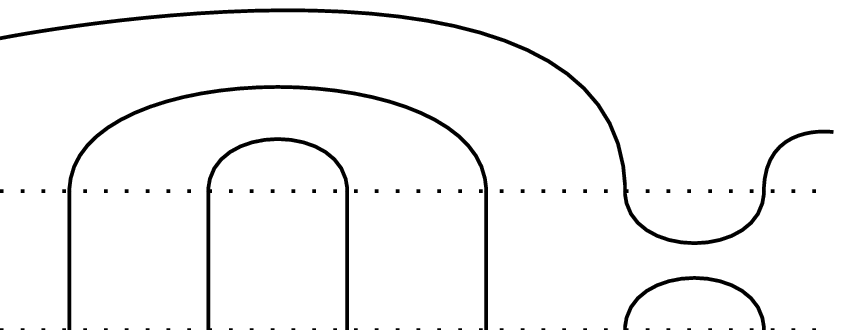}}
\put(130,10){$=$}
\put(150,10){\epsfxsize=120pt\epsfbox{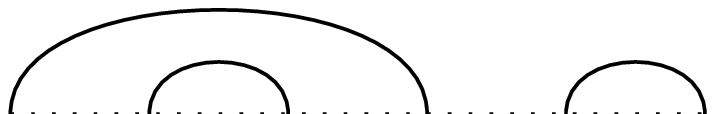}}
\end{picture}
\]
\end{example}

Instead of the Hamiltonian (\ref{eq:hamCdef}) which does not contain
boundary operators one may also consider the Hamiltonians $H^{\rm
M}_n: \Sp(\F^{\rm re}_{n}) \rightarrow \Sp(\F^{\rm re}_{n})$ and
$H^{\rm O}_n: \Sp(\F^{\rm e}_{n}) \rightarrow \Sp(\F^{\rm e}_{n})$ 
\begin{eqnarray*}
&&H^{\rm M}_n = (1-f_n) + \sum_{j=1}^{n-1} (1-e_j) \\
&&H^{\rm O}_n = (1-f_1) + (1-f_n) + \sum_{j=1}^{n-1} (1-e_j)
\end{eqnarray*} 
where M denotes {\bf mixed boundary conditions} \cite{MitraNGB02} and
O denotes {\bf open boundary conditions}. In the case of open
boundaries loop segments can ``penetrate the boundary'' and end on an
external point, as in Example \ref{ex:bound}. For closed boundaries
loops are not allowed to end on either external point, while for mixed
boundaries, the loops can only penetrate the right boundary but not
the left.  

\subsubsection{Periodic boundary conditions}
The action of the generators $e_j$ can be extended to directed
matchings using the graphical representation and keeping track of the
orderings of $i,j,k$ in (\ref{eq:TLgen_def}). Let us denote all
non-crossing directed (near-)perfect matchings of $[n]$ by
$\F_n^*$. For directed perfect matchings we can still use the
parentheses notation, e.g. $()\sim (1,2)$ and $)(\sim (2,1)$, and
graphically they can be conveniently depicted on a cylinder. A
matching $(i,j)$ for $i<j$ is represented by a loop segment over the
front of the cylinder and for $i>j$ by a loop segment over the
back. There is an equivalently graphical representation on an
annulus, the annulus being the top-view of the cylinder. The matching
$(i,j)$ is represented by a loop segment keeping the center of the
annulus to its left.
\begin{example}
\label{ex:dirmatch}
The directed matchings (1,2)(3,4)=()() and (1,2)(4,3)=())( are
graphically represented on the cylinder and annulus by
\[
\begin{picture}(210,120)
\put(0,90){$()()$}
\put(20,90){$\sim$}
\put(40,80){\epsfxsize=75pt\epsfbox{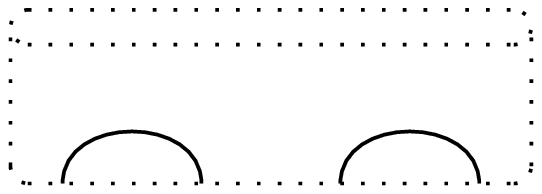}}
\put(46,70){$\sstyle 1$}
\put(65,70){$\sstyle 2$}
\put(86,70){$\sstyle 3$}
\put(105,70){$\sstyle 4$}
\put(130,90){$\sim$}
\put(150,70){\epsfxsize=50pt\epsfbox{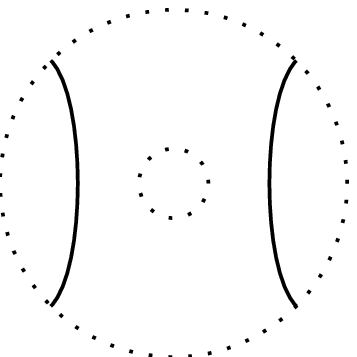}}
\put(150,115){$\sstyle 1$}
\put(150,70){$\sstyle 2$}
\put(195,70){$\sstyle 3$}
\put(195,115){$\sstyle 4$}
\put(0,20){$())($}
\put(20,20){$\sim$}
\put(40,10){\epsfxsize=75pt\epsfbox{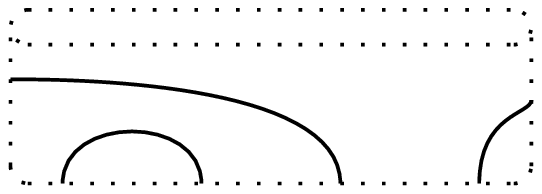}}
\put(46,0){$\sstyle 1$}
\put(65,0){$\sstyle 2$}
\put(86,0){$\sstyle 3$}
\put(105,0){$\sstyle 4$}
\put(130,20){$\sim$}
\put(150,0){\epsfxsize=50pt\epsfbox{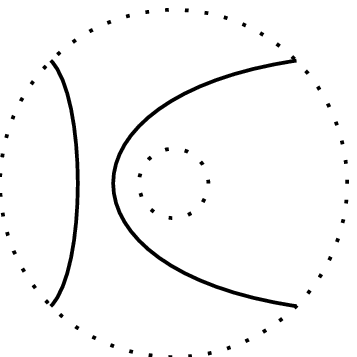}}
\put(150,45){$\sstyle 1$}
\put(150,0){$\sstyle 2$}
\put(195,0){$\sstyle 3$}
\put(195,45){$\sstyle 4$}
\end{picture}
\]
\end{example}
The definitions (\ref{eq:TLgen_def}) for the generators $e_j$
extended to directed matchings are,
\begin{equation}
e_j: \left\{ 
\begin{array}{lcl}
(j,j+1),\; (j+1,j) & \mapsto & (j,j+1) \\
(i,j)(j+1,k) & \mapsto & (i,j)(j,j+1) \\
(j,i)(j+1,k),\; (i,j)(k,j+1) & \mapsto & (k,i)(j,j+1)
\end{array}\right.,
\label{eq:TLgen_perdefj}
\end{equation}
for $i<j<k$.

If $n$ is odd, there will be one unpaired vertex in
$\F^*_n$. In the parentheses notation we denote this by a vertical
line, e.g. $(1,2)\{3\} \sim ()|$, while graphically we think of it as
a defect line running from the top to the bottom of the cylinder. The
action of $e_j$ extended to near-perfect directed matchings is
\begin{equation}
e_j:\;\; (i,j)\{j+1\},\; \{j\}(j+1,i) \;\; \mapsto \;\; \{i\}(j,j+1). 
\label{eq:TLgen_perdefj2}
\end{equation}

Graphically it is natural to introduce an additional generator $e_n$
for directed matchings, see Levy \cite{Levy91} and Martin
\cite{Martin91}. Its non-trivial action on elements $F\in\F_n^*$
containing $1$ or $n$ is similar to (\ref{eq:TLgen_perdefj}) and is
described by,
\begin{equation}
e_n: \left\{ 
\begin{array}{lcl}
 (n,1),\; (1,n) & \mapsto & (n,1) \\
(1,i)(j,n) & \mapsto & (j,i)(n,1) \\
(1,i)(n,j),\; (i,1)(j,n) & \mapsto & (i,j)(n,1) \\
\{1\}(i,n),\; (1,i)\{n\} & \mapsto & \{i\}(n,1)
\end{array}\right.
\label{eq:TLgen_perdefn}
\end{equation}
\begin{lemma}
The generator $e_n$ satisfies the relations
\begin{eqnarray*} &&e_n^2 = e_n,\qquad e_n e_1 e_n = 1,\qquad e_1 e_n e_1 =
e_1 \\
&&e_ne_j = e_je_n \qquad j\notin\{1,n-1\},
\end{eqnarray*}
and the generators $e_j,\;j\in\{1,\ldots,n\}$ can be represented
graphically on a cylinder as
\[
e_j \quad =\quad
\begin{picture}(140,15)
\put(0,-10){\epsfxsize=140pt\epsfbox{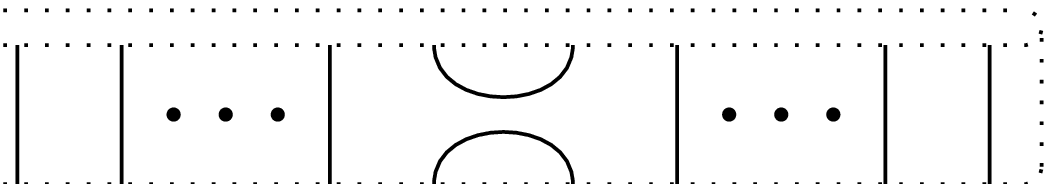}}
\put(.5,-18){$\sstyle 1$}
\put(15,-18){$\sstyle 2$}
\put(38,-18){$\sstyle j-1$}
\put(58,-18){$\sstyle j$}
\put(72,-18){$\sstyle j+1$}
\put(89,-18){$\sstyle j+2$}
\put(117,-18){$\sstyle n-1$}
\put(135.5,-18){$\sstyle n$}
\end{picture}\hspace{10pt}
\]
\vskip20pt
\noindent
\end{lemma}

The definitions (\ref{eq:TLgen_perdefj}) and (\ref{eq:TLgen_perdefn})
imply that for perfect matchings non-contractible loops running around
the cylinder can be erased. For near-perfect matchings they allow to
neglect the winding of the defefct line. A thorough discussion how to
capture this in algebraical rules is given in \cite{MartinS93}, see
also \cite{PearceRGN02}.
\begin{example} 
$e_3 ())( \; = ()() \in \F^*_4$ is displayed graphically as 
\[
\begin{picture}(185,45)
\put(0,0){\epsfxsize=75pt\epsfbox{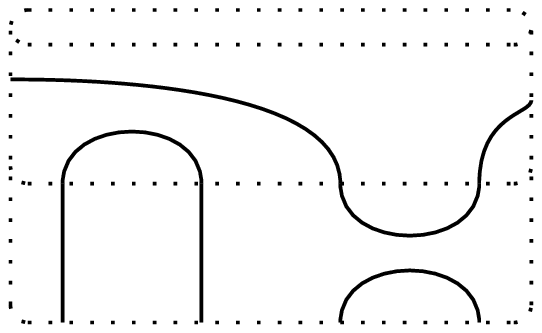}}
\put(90,10){$=$}
\put(110,0){\epsfxsize=75pt\epsfbox{wordcyl1.eps}}
\end{picture}
\]
and $e_3 ()| \; =\; )|(\; \in \F^*_3$ corresponds to
\[
\begin{picture}(185,50)
\put(0,0){\epsfxsize=60pt\epsfbox{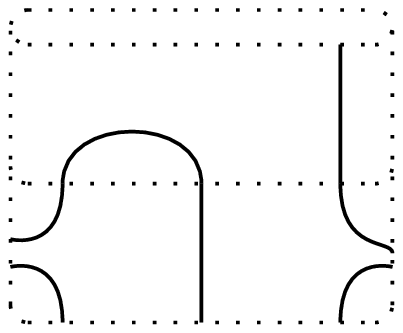}}
\put(90,10){$=$}
\put(110,0){\epsfxsize=30pt\epsfbox{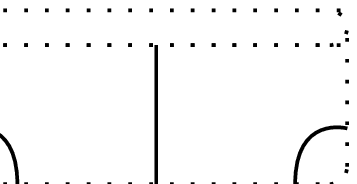}}
\end{picture}
\]
\end{example}

We define a new Hamiltonian $H_n^{\rm P*}: \Sp(\F^*_n) \rightarrow \Sp(\F^*_n)$ by
\begin{equation}
H_n^{\rm P*} = \sum_{j=1}^{n}(1-e_j), \label{eq:hamP*def}
\end{equation}
where P stands for {\bf periodic boundary conditions}. For even $n$,
the Hamiltonian (\ref{eq:hamP*def}) also has an action on non-directed
non-crossing perfect matchings. Graphically this corresponds to
closing the top of the cylinder, or removing the inner disk of the
annulus. The two distinct matchings in Example \ref{ex:dirmatch} then
become equal. This case we denote by $H_{2n}^{\rm P}: \Sp(\F_{2n})
\rightarrow \Sp(\F_{2n})$ 
\[
H_{2n}^{\rm P} = \sum_{j=1}^{2n}(1-e_j). 
\]

\section{Fully packed loop diagrams}
\label{se:FPL}
The {\bf fully packed loop model} (FPL) \cite{BatchBNY96,Wiel00} is a
model of polygons on a lattice such that each vertex is visited
exactly once by a polygon. The model is an alternative representation of
the {\bf six-vertex model} on the square lattice, as will be described
below. Here we consider the six-vertex model with {\bf domain wall}
boundary conditions which were introduced by Korepin \cite{Kore82}, as
well as related boundary conditions.   

It is well known that ASMs are in bijection with configurations of the
six-vertex model with domain wall boundary conditions, see for example
Elkies et al. \cite{ElkiesKLP92}. The correspondence between entries
in an ASM and the six vertex configurations is given by 
\[
\begin{picture}(200,60)
\put(0,20){\epsfxsize=200pt\epsfbox{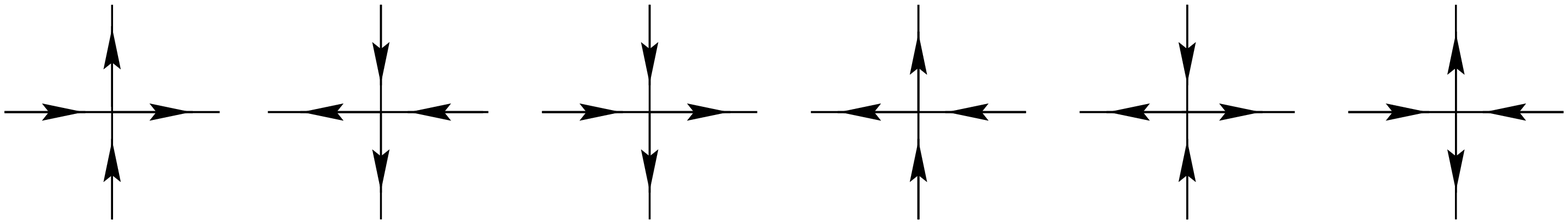}}
\put(12,0){$0$}\put(46,0){$0$}\put(80,0){$0$}
\put(115,0){$0$}\put(145,0){$-1$}\put(184,0){$1$}
\end{picture}
\]
\begin{example}
There are seven $3\times3$ ASMs. Six of these are the $3\times3$
permutation matrices and only one of them contains a $-1$. Its
corresponding six-vertex configuration is given by
\[
\begin{picture}(170,60)
\put(0,26){$
\left(\begin{array}{@{}ccc@{}}
0 & 1 & 0 \\
1 & -1 & 1 \\
0 & 1 & 0
\end{array}\right)$}
\put(80,26){$\leftrightarrow$}
\put(110,0){\epsfxsize=60pt\epsfbox{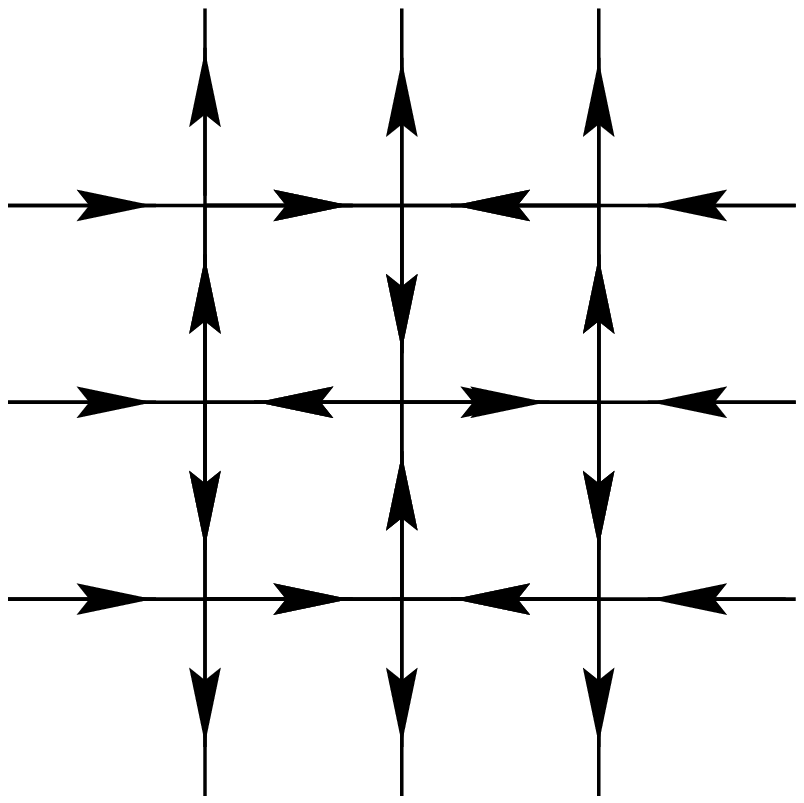}}
\end{picture}
\]
\end{example}

The six-vertex model was used extensively by Kuperberg \cite{Kupe96,Kupe00}
to prove many enumerations of ASMs and their
symmetry classes using and generalizing the Izergin and Tsuchiya
determinants \cite{Izergin82,IzerCK92,KoreIB93,Tsuch98}. It was
also used by Zeilberger \cite{Zeilb96} in his proof of the refined
alternating-sign matrix conjecture.

A fully packed loop configuration is obtained from a six-vertex
configuration by dividing the square lattice into its even and odd
sublattice denoted by $A$ and $B$ respectively. Instead of arrows, only those
edges are drawn that on sublattice $A$ point inward and on sublattice $B$
point outward.
\[
\begin{picture}(230,106)
\put(30,0){\epsfxsize=200pt\epsfbox{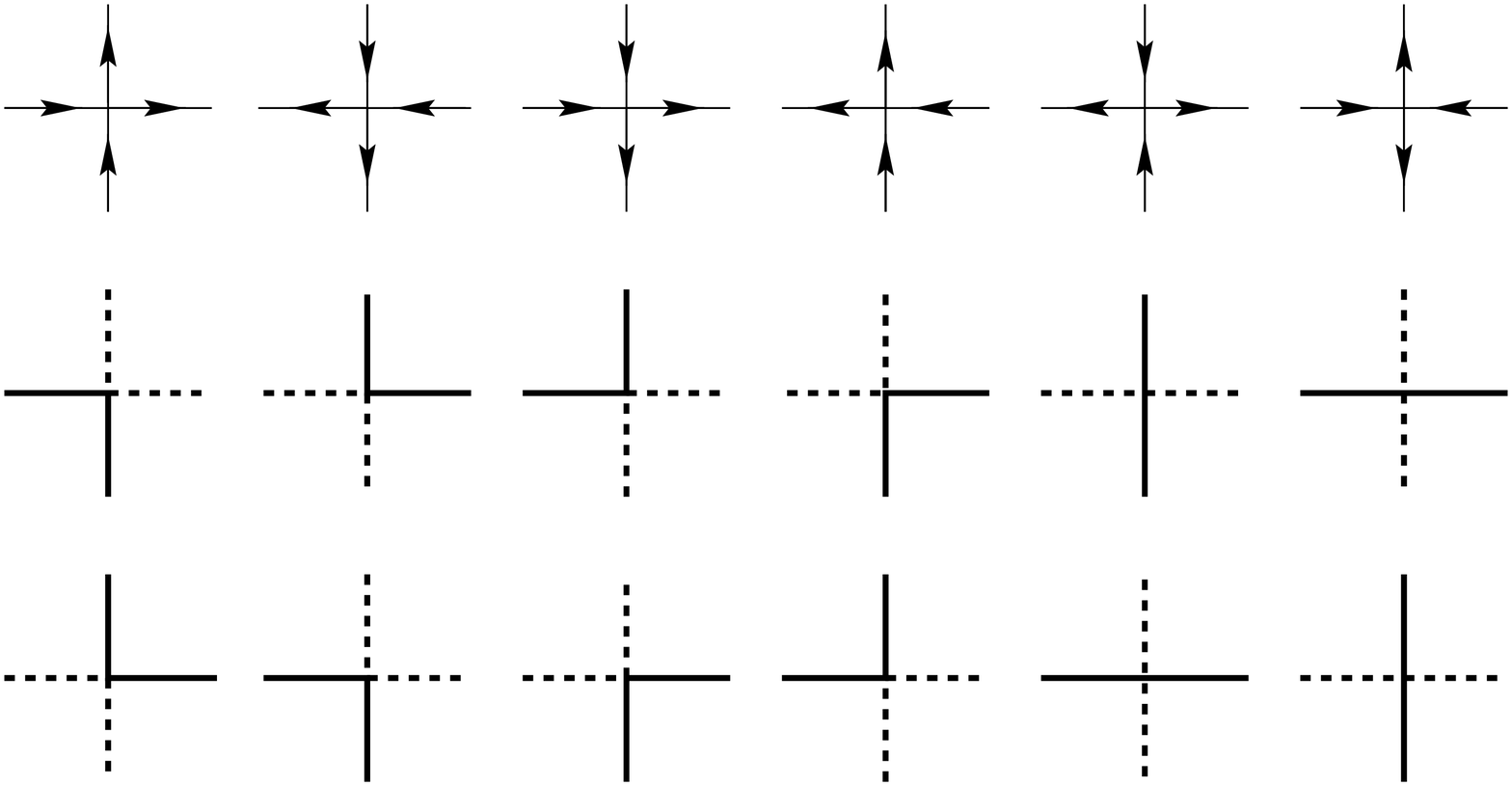}}
\put(0,50){$A$}
\put(0,12){$B$}
\end{picture}
\]
We take the vertex in the upper left corner to belong to sublattice $A$.
The domain wall boundary condition translates into a boundary
condition for the loops. Loops either form closed circuits, or begin
and end on boundary sites which are prescribed by the boundary in- and
out-arrows on sublattice $A$ and $B$ respectively. Such sites are
called {\bf designated boundary sites}. From the above it
follows that ASMs are in bijection with fully packed loop diagrams
on a square grid. The {\bf grid} for $n\times n$ 
ASMs is denoted by $G_n$. 
\begin{example}
The $7\times7$ ASMs are in bijection with FPL diagrams on $G_7$
\vskip2mm
\[
\begin{picture}(140,100)
\put(0,50){$G_7 =$}
\put(40,0){\epsfxsize=100pt\epsfbox{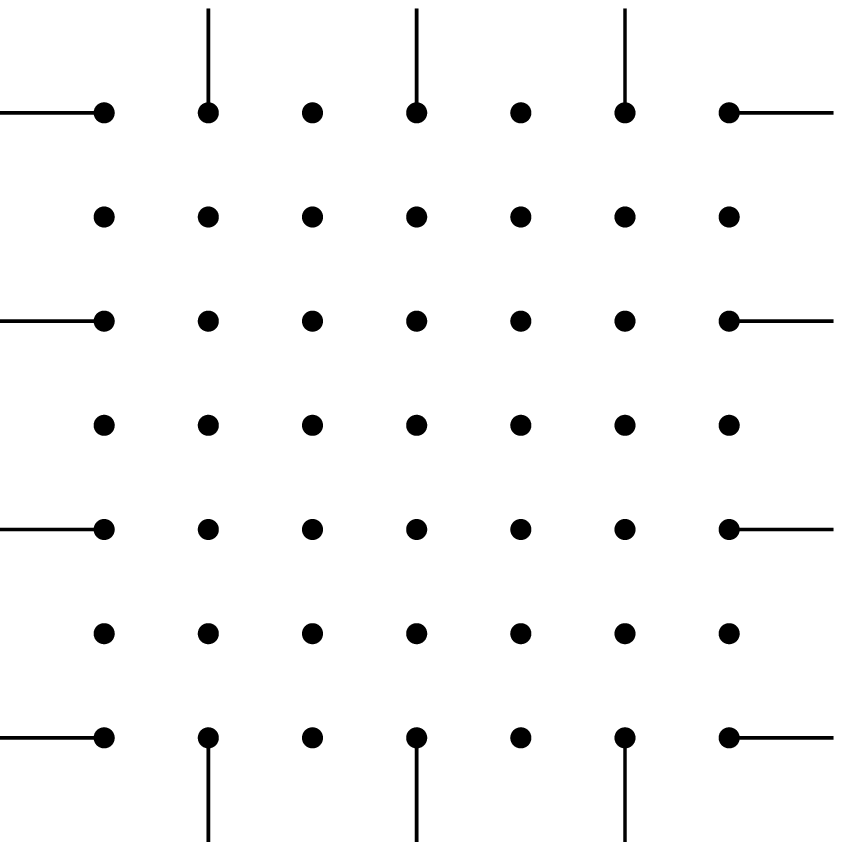}}
\end{picture}
\]
\end{example}

\subsection{Matchings of FPL diagrams}
In each FPL diagram every designated boundary site is connected
to another such site. Number the designated boundary sites, starting with
$1$ at the left boundary in the upper left corner and ending with $2n$
at the top boundary in the same corner. We then have
\begin{lemma}
Each FPL diagram on a grid $G_n$ defines a perfect matching of
$[2n]$. A matching therefore defines an equivalence relation on FPL
diagrams. 
\end{lemma}

The {\bf ($G$,$F$)-cardinality} $M_F(G)$ is the number of FPL
diagrams with matching $F$ on grid $G$.
\begin{example}
There are $7$ FPL diagrams on $G_3$,
\[
\begin{picture}(300,111)
\put(0,0){\epsfxsize=300pt\epsfbox{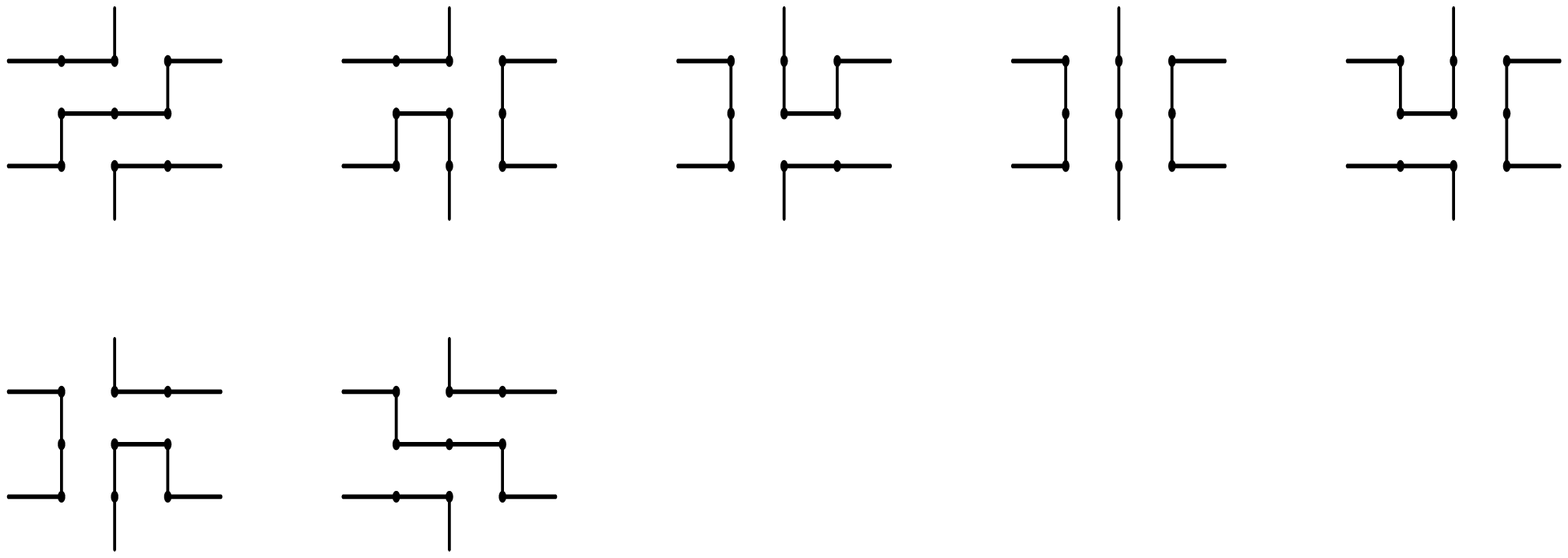}}
\end{picture}
\]
The $2$nd and $5$th diagram have the same matching and likewise with
the $3$rd and $6$th diagram. The others have distinct matchings. The
set of matchings is therefore $\{)(())(,)()()(,()()(),()(()),(())()\}$
and the corresponding set of cardinalities is $\{1,2,2,1,1\}$.
\end{example}

\subsection{Symmetry classes}
Requiring ASMs to be symmetric with respect to one of the symmetries
of the square puts constraints on the loop configurations at certain
edges. For $(2n+1)\times(2n+1)$ {\bf vertically symmetric} ASMs there 
must be a loop running from top to bottom along the symmetry axis, and
it is also not hard to see that the edges in the top and bottom row
are fixed. Requiring vertical symmetry reduces the configuration space
of FPL diagrams on $G_{2n+1}$ to that of FPL diagrams on a rectangular
grid of size $n\times (2n-1)$. Let $G_{2n}^{\rm V}$ denote this grid. 
\begin{lemma}
Each FPL diagram on a grid $G_{2n}^{\rm V}$ defines a perfect matching on $[2n]$.
\end{lemma}

\begin{example} 
The $26$ vertically symmetric ASMs of size $7\times 7$ are in
bijection with FPL diagrams on $G_{6}^{\rm V}$,
\[
\begin{picture}(260,100)
\put(0,40){$G^{\rm V}_6=$}
\put(35,20){\epsfxsize=90pt\epsfbox{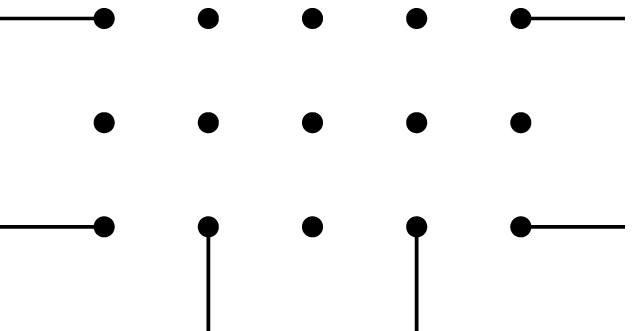}}
\put(135,48){$\leftarrow$}
\put(160,0){\epsfxsize=100pt\epsfbox{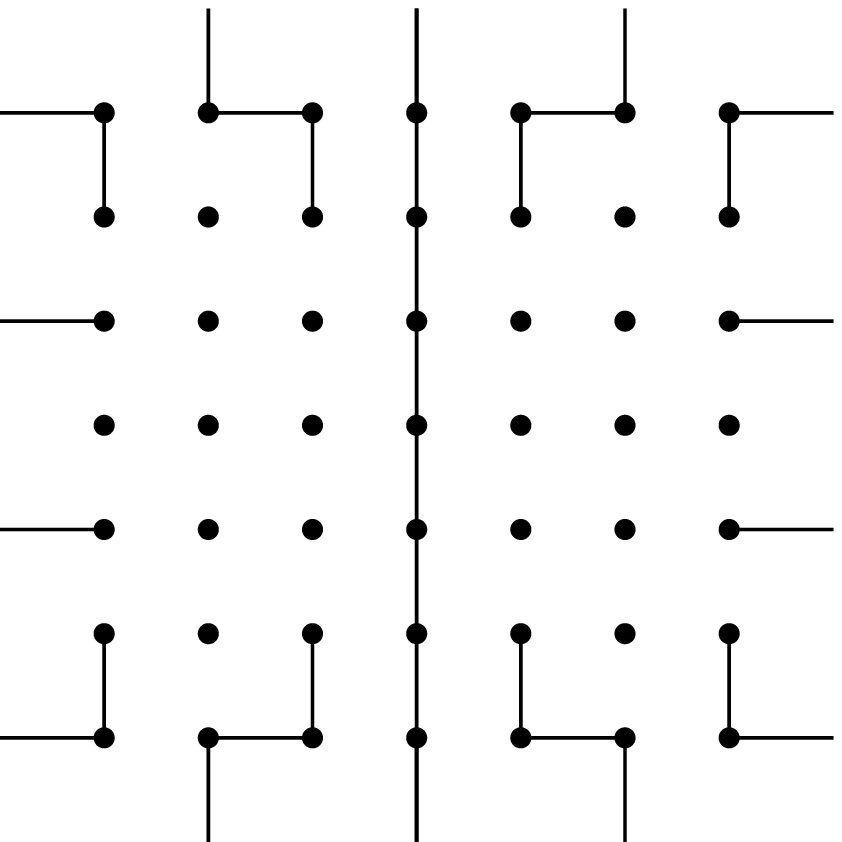}}
\end{picture}
\]
\end{example}

A class of FPL diagrams closely related to those on $G_{2n}^{\rm V}$
are diagrams with an odd number of designated boundary sites.
This class is defined by FPL diagrams on an $n\times (2n+1)$ grid,
denoted by $G^{\rm V}_{2n+1}$, such that the unpaired loop line may
end anywhere on the top boundary. 

\begin{lemma}
Each such FPL diagram on a grid $G_{2n+1}^{\rm
V}$ defines a near-perfect matching on $[2n+1]$, or equivalently, a
perfect right extended $(n,1)$-matching on $[2n+1]$.
\end{lemma}

\begin{remark}
The upper boundary of $G_{2n+1}^{\rm V}$ plays a role analogous to the
vertex $\{n+1\}$.
\end{remark}

\begin{example}
An FPL diagram on $G^{\rm V}_7$ with matching $()()()($
is 
\[
\begin{picture}(120,50)
\put(0,0){\epsfxsize=120pt\epsfbox{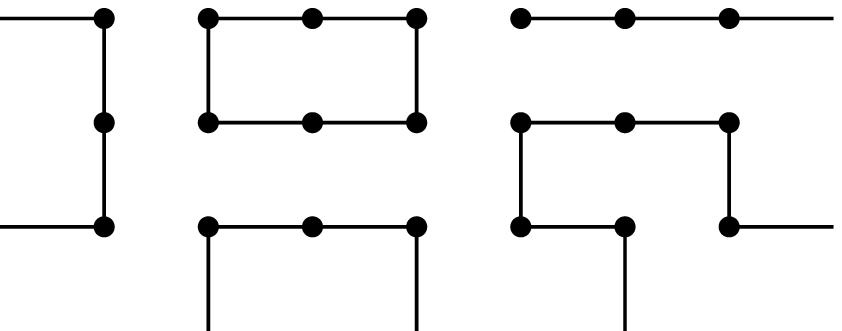}}
\end{picture}
\]
\end{example}

For $(2n+1)\times(2n+1)$ {\bf vertically and horizontally symmetric} ASMs
the horizontal symmetry axis contains loop segments as well as all
horizontal edges crossing the vertical symmetry axis. Furthermore,
the complete boundary layer is fixed. Requiring horizontal and
vertical symmetry reduces the configuration space of FPL diagrams on
$G_{2n+3}$ to that of FPL diagrams on a square grid of size
$n\times n$, which is denoted by $G_{n}^{\rm VH}$.
\begin{lemma}
Each FPL diagram on a grid $G_{n}^{\rm VH}$ defines a perfect
right extended matching on $[n]$.
\end{lemma}

\begin{remark}
The boundary sites on $G_{n}^{\rm VH}$ coming from the horizontal
edges on the vertical symmetry axis of $G_{2n+3}$ play the role of the
extra site at $\{n+1\}$.
\end{remark}

\begin{example} 
The $6$ vertically and horizontally symmetric ASMs of size $9\times 9$
are in bijection with FPL diagrams on $G_{3}^{\rm VH}$
\[
\begin{picture}(260,130)
\put(0,55){$G^{\rm VH}_3=$}
\put(35,30){\epsfxsize=60pt\epsfbox{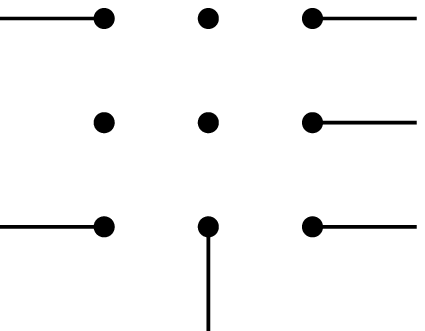}}
\put(110,55){$\leftarrow$}
\put(140,0){\epsfxsize=120pt\epsfbox{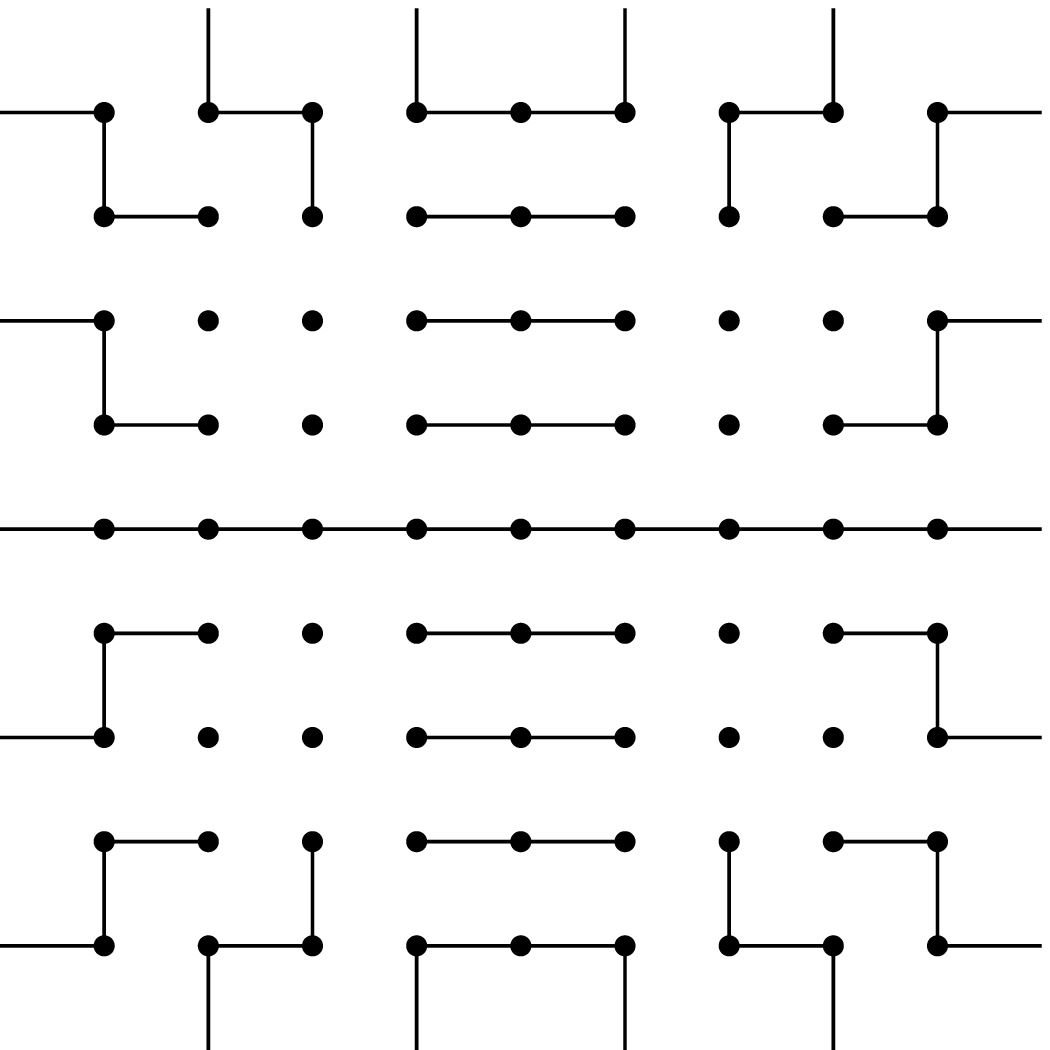}}
\end{picture}
\]
\end{example}

{\bf Half-turn symmetric} ASMs, or HTASMs, are symmetric under
rotation by $180$ degrees. For HTASMs of size $2n$ one has only to consider
FPL diagrams on the lower half. In contrast to the case of vertically
symmetric ASMs, the sites on the top boundary of this rectangular grid
can now also be connected via one of $n $ arcs, or HT boundaries 
\cite{Kupe00}. Loops crossing the horizontal symmetry axis of the
square, map to loops on such arcs. Requiring half turn symmetry reduces
the configuration space of FPL diagrams on $G_{2n}$ to that of FPL
diagrams on a rectangular grid of size $2n^2$ with HT boundaries, which is denoted by
$G_{2n}^{\rm HT}$,  
\[
\begin{picture}(270,90)
\put(0,40){$G^{\rm HT}_6=$}
\put(40,20){\epsfxsize=90pt\epsfbox{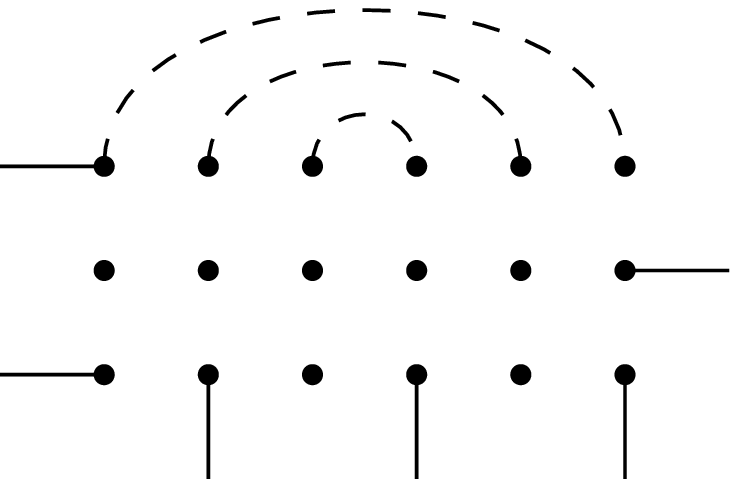}}
\put(150,40){$\leftarrow$}
\put(180,0){\epsfxsize=90pt\epsfbox{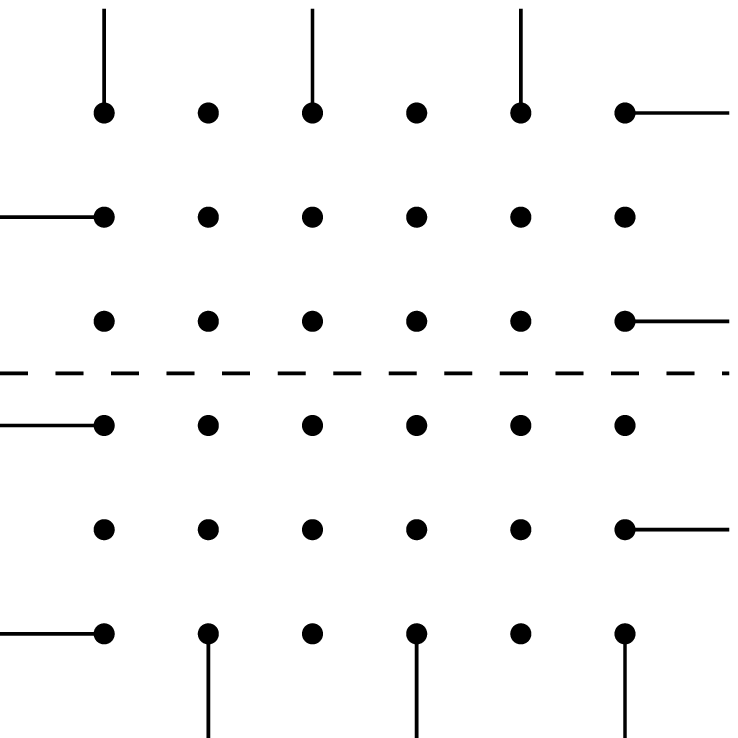}}
\end{picture}
\]
\begin{example}
An FPL diagram on $G^{\rm HT}_6$ with matching $)()()($
is 
\[
\begin{picture}(90,50)
\put(0,0){\epsfxsize=90pt\epsfbox{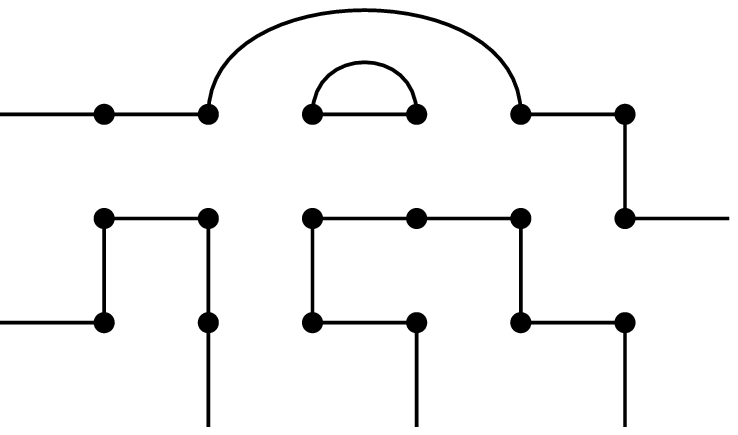}}
\end{picture}
\]
\end{example}
For odd-sized HTASMs the situation is a little more complicated. In
this case there is one loop segment running across the FPL diagram
dividing it in two identical parts. However, the shape of these parts
depend on the way the loop segment runs across the cylinder and
there is no unique reduced graph to accommodat for all possible link
patterns. In the is case we define $G_{2n+1}^{\rm HT}$ as $G_{2n+1}$,
the graph for odd unrestricted FPL diagrams, with the proviso that it
has to be half-turn symmetric. 
\begin{example}
An FPL diagram on $G^{\rm HT}_5$ with matching $)()|($
is 
\[
\begin{picture}(90,90)
\put(0,0){\epsfxsize=90pt\epsfbox{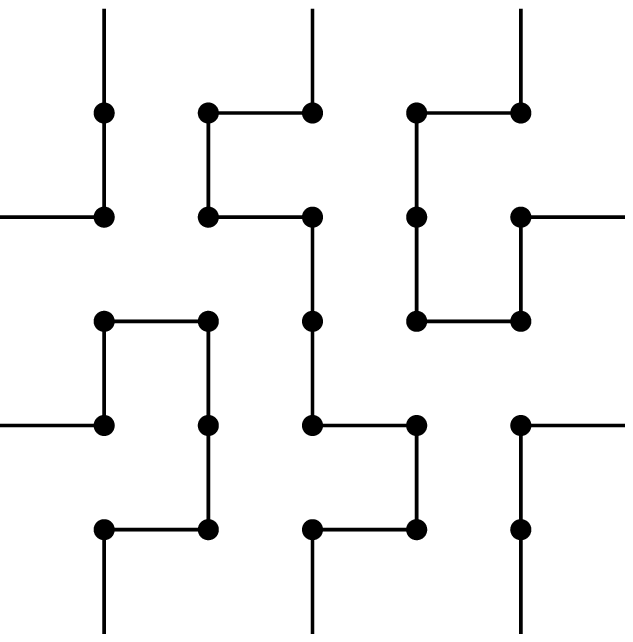}}
\end{picture}
\]
\end{example}

\subsection{Boundary conditions and symmetry classes}
Finally we can state a conjecture that relates the stationary state
for different boundary conditions in the loop model to symmetry
classes of FPL diagrams. This conjecture is a collection of results
obtained for special cases by several authors
\cite{BatchGN01,MitraNGB02,PearceRGN02,RazuS01a,RazuS01b,RazuS01c}.
\begin{conj}
\label{co:TL-FPL}
Let $\{H,\F,G\}$ be any of the following four triples,
\[
\renewcommand{\arraystretch}{1.5}
\begin{array}{c|c|c}
H & \F & G \\\hline
H^{\rm C}_{n} & \F_{n} & G^{\rm V}_{n}\\
H^{\rm M}_{n} & \F^{\rm re}_{n} & G^{\rm VH}_n\\
H^{\rm P}_{2n} & \F_{2n} & G_n\\
H^{\rm P*}_{n} & \F^*_{n} & G^{\rm HT}_{n}\\
\end{array}
\]
and let $P \in\Sp(\F)$ be a solution of 
\[
H P = 0,\qquad P=\sum_{F\in \F} a_F F,
\]
then $a_F$ is equal to the $(F,G)$-cardinality $M_F(G)$, i.e., the
number of fully packed loop diagrams on $G$ with matching $F\in\F$. 
\end{conj} 
There are also some conjectures for open boundary conditions
\cite{MitraNGB02}, but it is not known to which grid that case
corresponds.
\begin{example}
The stationary state for closed boundaries and $n=6$ is given by $P_6^T
= (11,5,5,4,1)$ on the basis $\F_6 =
\{()()(),(())(),()(()),(()()),((()))\}$, see Example
\ref{ex:Cn=6}. Indeed one finds among the $26$ FPL diagrams on $G^{\rm
V}_6$ precisely $11$ with matching $()()()$, $5$ with matching
$(())()$, $5$ with matching $()(())$, $4$ with matching $(()())$ and
$1$ with matching $((()))$. The $11$ diagrams with matching $()()()$
are printed bold, 
\vskip2mm
\[
\begin{picture}(340,200)
\put(0,0){\epsfxsize=340pt\epsfbox{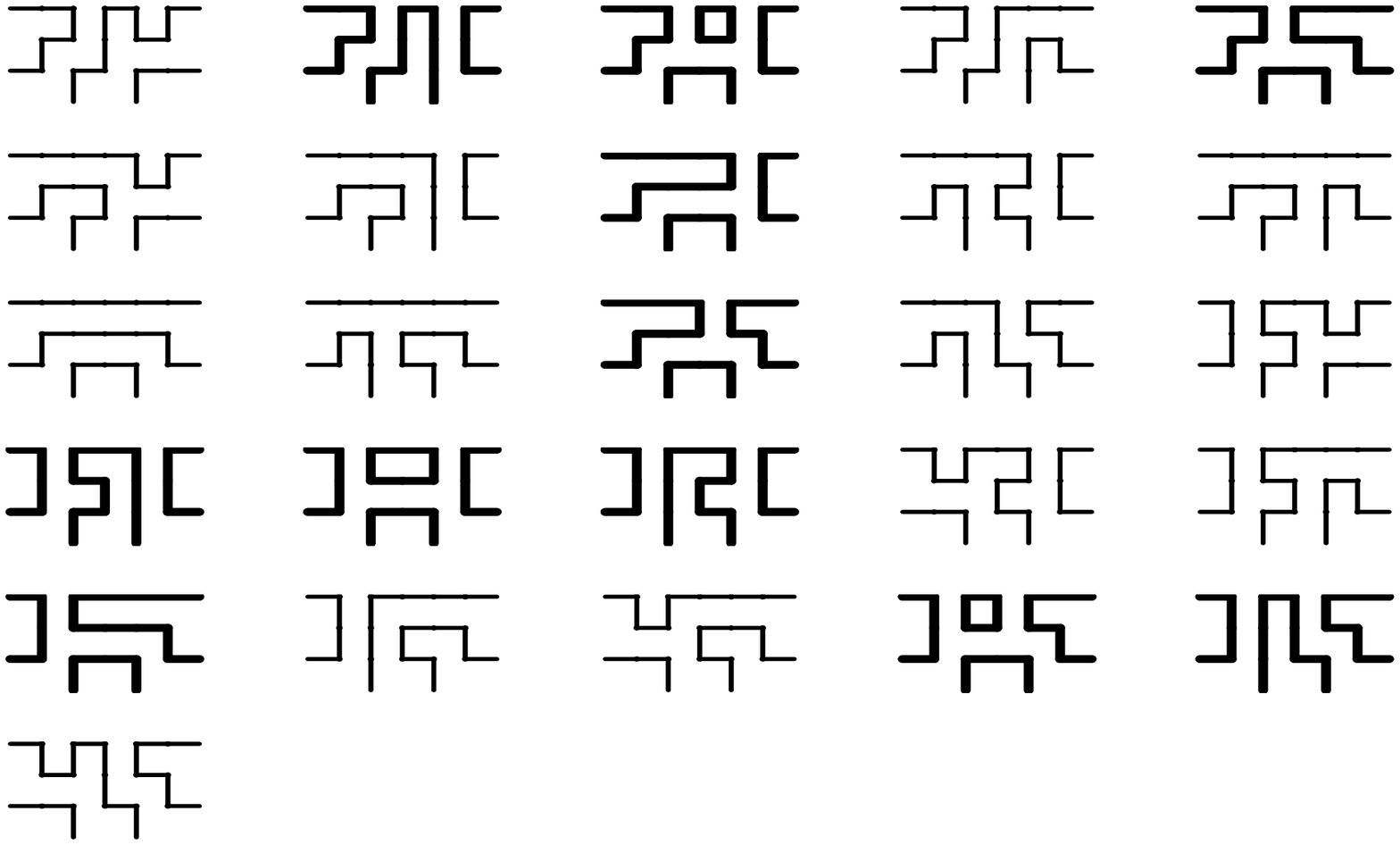}}
\end{picture}
\]
\end{example}

\section{Correlation functions and expectation values}
\label{se:Nests}
In physics one would rather like to know expectation
values and correlation functions than
the stationary state itself. An expectation value is nothing else but
the ratio of a refined to a full enumeration. If the refinement is
with respect to two or more constraints one usually speaks of a
correlation function rather than expectation value.
A well known example of an expectation value is given by the refined ASM conjecture,
which was formulated by Mills et al.
\cite{MillsRR82,MillsRR83} and proved by Zeilberger
\cite{Zeilb96}. Recent progress has been made by Stroganov for 
double refined ASM correlations \cite{Strog02}. 

For other correlation functions or expectation values 
closed formulae have been conjectured on the basis of explicit
calculations for small system sizes. The formula for the 
emptiness formation probability (see e.g. \cite{KoreIB93}) in the XXZ spin 
chain representation for periodic and twisted boundary conditions was
conjectured by Razumov and Stroganov in \cite{RazuS01a} and
\cite{RazuS01d}. The asymptotic behaviour of that conjecture was
proved by Kitanine et al. in \cite{KitaMST02}
using the  integrability of the XXZ spin chain. In the loop representation
several other conjectures for expectation values are stated in
\cite{MitraNGB02,Zuber03}. 

Below we will give examples of explicit formulas for
nest distribution functions in the case of closed and periodic
boundary conditions, or vertically and half-turn symmetric FPL diagrams
respectively.

\subsection{Nests}
The {\bf nests} of a non-crossing perfect matching $F$ are the
sequences of consecutive parentheses of $F$ that close among
themselves and are not enclosed by other parentheses. The number of nests of 
$F$ is denoted by $c_F$. In terms of Dyck paths (see Remark
\ref{rem:Dyck}), the number of nests is one less the number of times
the path touches the real axis.  
\begin{example}
The number of nests for the non-crossing perfect matchings $F =
()()(),(())(),()(()),(()()),((()))$ is $c_{F}=3,2,2,1,1$.
\end{example}
The nests of an FPL diagram are the nests defined by its
matching.

Let $\F_G$ denote the set of non-crossing perfect matchings defined by
the FPL diagrams on the grid $G$. The {\bf nest distribution function}
$P_G(k)$ is defined as the number of FPL diagrams on $G$ with $k$
nests, 
\[
P_{G}(k) = \sum_{F \in \F_G} M_F(G) \delta_{k,c_F},\qquad Z_G = \sum_k
P_G(k),
\]
and the {\bf average number of nests} is defined as 
\be
\langle k \rangle_G = \frac{1}{Z_G} \sum_{k} k\, P_{G}(k).
\label{eq:cav}
\ee
Below we show that in some cases an exact formula for these quantities
can be guessed from explicit data for small sizes. As a rule,
a number can only be guessed easily if it factorises into small
primes. It should be mentioned that in such cases there are several
helpful software utilities that could be used, see e.g. Appendix A in
\cite{Kratt99}.
\begin{ta}
\label{ta:HTcav}
The average number of nests on $G=G^{\rm HT}_{2n}$ for $n \in
\{2,\ldots,7\}$ and the prime factorisations of its
denominator and its numerator.
\[
\renewcommand{\arraystretch}{1.6}
\begin{array}{c|c|c|c|c|c|c}
2n & 4 & 6 & 8 & 10 & 12 & 14 \\ \hline
\langle k\rangle_{G} & \frac{8}{5} & \frac{21}{10} & \frac{28}{11} &
\frac{65}{22} & \frac{624}{187} & \frac{3485}{935}\\ \hline
{\rm den} & 5 & 2\cdot 5 & 11 & 2\cdot 11 & 11\cdot 17 & 5\cdot
11\cdot 17 \\ \hline
{\rm num} & 2^3 & 3\cdot 7 & 2^2\cdot 7 & 5\cdot 13 & 2^4\cdot 3\cdot
13 & 2\cdot 7\cdot 13\cdot 19
\end{array}
\]
\end{ta}
It is not difficult to derive a formula for $\langle
k\rangle_{G^{\rm HT}_{2n}}$ that exactly reproduces the numbers in 
Table \ref{ta:HTcav}. We have
\begin{conj}
\label{co:HTcav}
The average number of nests on $G=G^{\rm HT}_{2n}$ is 
\[
\langle
k\rangle_{G} = n \prod_{j=1}^{n-1} \frac{3j+1}{3j+2}
\sim \frac{ \Gamma(5/6)}{\sqrt{\pi}}\; (2n)^{2/3}\qquad
(n\rightarrow \infty). 
\]
\end{conj}
This guessing does not always work as we can see from the next table
where the same quantities are given for $G^{\rm V}_{2n}$,
\begin{ta}
\label{ta:Vcav}
The average number of nests on $G=G^{\rm V}_{2n}$ for $n \in
\{2,\ldots,7\}$ and the prime factorisations of its
denominator and its numerator.
\[
\renewcommand{\arraystretch}{1.6}
\begin{array}{c|c|c|c|c|c|c}
2n & 4 & 6 & 8 & 10 & 12 & 14 \\ \hline
\langle k\rangle_{G} & \frac{5}{3} & \frac{29}{13} & \frac{52}{19} &
\frac{913}{285} & \frac{1693}{465} & \frac{69769}{17205}\\ \hline
{\rm den} & 3 & 13 & 19 & 3\cdot 5\cdot 19 & 3\cdot 5\cdot 31 & 3\cdot
5\cdot 31\cdot 37 \\ \hline
{\rm num} & 5 & 29 & 2^2\cdot 13 & 11\cdot 83 & 1693 & 7\cdot 9967
\end{array}
\]
\end{ta}
The appearance of large prime factors in the numerator makes it hard
to conjecture a formula. However, for the case of nests we are very
fortunate because while we may not be able to guess the
average (\ref{eq:cav}) for $G=G^{\rm V}_{2n}$, we can find the
{\em complete} distribution function.

\begin{conj}
\label{co:Vnestdf}
Let $(a)_k = \Gamma(a+k)/\Gamma(a)$ and $G=G^{\rm V}_{2n}$. The nest
distribution function $P_G(k)$ is given by
\[
P_G(k) = k \frac{4^{n+k}}{27^n} \frac{(1/2)_{n+k}}{(1/3)_{2n}}
\frac{(3n+1)!(2n-k-1)!}{n!(n-k)!(2n+k+1)!}
A^{\rm V}_{2n+1}, 
\]
where
\[
A^{\rm V}_{2n+1} = \prod_{j=0}^{n-1} (3j+2)
\frac{(6j+3)!(2j+1)!}{(4j+2)!(4j+3)!} \; =\; 1,3,26,646,\ldots
\]
is the number of $(2n+1) \times (2n+1)$ vertically symmetric ASMs.
\end{conj}
This conjecture has been checked for $n$ up to $8$.
Assuming the conjecture the average number of nests can be
calculated. It turns out that (\ref{eq:cav}) with 
Conjecture \ref{co:Vnestdf} is summable (see below) and we find,
\begin{cor}
\label{cor:Vcav}
The average number of nests on $G=G^{\rm V}_{2n}$ is
\begin{eqnarray*}
\langle k\rangle_{G} &=& \frac{1}{A^{\rm V}_{2n+1}}\sum_{k=1}^n
k P_G(k) = \frac{1}{3}\left(
\prod_{j=0}^{n-1}\frac{(2j+1)(3j+4)}{(j+1)(6j+1)} -1 \right)\\
&\sim& \frac{\Gamma(1/3)
\sqrt{3}}{2 \pi} (2n)^{2/3}\qquad (n\rightarrow \infty). 
\end{eqnarray*}
\end{cor}
Corollary \ref{cor:Vcav} fits the data in Table \ref{ta:Vcav}. While we
obtain a nice formula for $\langle k\rangle_{G^{\rm V}_{2n}}$, the
simple subtraction of $1/3$ in Conjecture \ref{cor:Vcav}) gives rise to
large prime factors in Table \ref{ta:Vcav} and makes it virtually
impossible to guess the result directly from there. 
A result similar to Conjecture \ref{co:Vnestdf} is obtained for nests
on $G^{\rm HT}_{2n}$.
\begin{conj}
\label{co:HTnestdf}
Let $G=G^{\rm HT}_{2n}$. The nest distribution function $P_G(k)$ is given by
\begin{equation}
P_G(k) = 3nk \frac{(2k)! (n+k-1)! (2n-k-1)!}
{k!^2(n-k)!(2n+k)!}A^{\rm HT}_{2n}.
\end{equation}
where
\[
A^{\rm HT}_{2n} = A_n^2 \prod_{j=0}^{n-1} \frac{3j+2}{3j+1} = 2,10,140,\ldots.
\]
is the number of $2n \times 2n$ half-turn symmetric ASMs. $A_n$ is
the number of $n\times n$ ASMs defined in
(\ref{eq:ASMno}).
\end{conj}
This conjecture has again been checked for $n$ up to $8$ and is consistent
with Conjecture \ref{co:HTcav}.

\subsection{Strange evaluations}
Corollary \ref{cor:Vcav} has been obtained from Conjecture \ref{co:Vnestdf}
using the Mathematica implementation of the Gosper-Zeilberger
algorithm \cite{Gosper78,Zeilb90,Zeilb91} by Paule and Schorn
\cite{PauleS95}. Using a similar procedure, the consistency of
Conjectures \ref{co:HTcav} and \ref{co:HTnestdf} was checked. Both
results can also be derived using the following hypergeometric
summation formula,\footnote{We are grateful to Christian Krattenthaler for
pointing this out.}
\be
{}_5F_4 \left[ 
\begin{array}{@{}c}
a,1+2a/3,1-2d,1/2+a+m,-m \\
2a/3,1/2+a+d,-2m,1+2a+2m 
\end{array};4\right] = \frac{(1-d)_m (a+1)_m} {(1/2)_m (1/2+a+d)_m},
\label{eq:strange}
\ee
which can be derived from one of the strange evaluations of Gessel and
Stanton \cite{GesselS82}. A detailed derivation is beyond the scope of
this paper, but one may for example derive Conjecture \ref{co:HTcav}
by writing (\ref{eq:cav}) in hypergeometric notation using Conjecture
\ref{co:HTnestdf}, and taking $a=3/2$, $m=n-1$ and $d=-1/3$ in 
(\ref{eq:strange}). The derivation of Corollary \ref{cor:Vcav} is
similar but a bit more complicated since one has to use an additional
contiguous relation.

\section{Hexagons with cut off corners}
\label{se:Hex}
The following conjecture stated by Mitra et al. \cite{MitraNGB02} claims
that certain elements of the stationary state of $H^{\rm C}_{2n}$ are
given by enumerations of hexagons with cut off corners calculated by
Ciucu and Krattenthaler \cite{CiucuK02}. 
Using the notation of \cite{MitraNGB02}, we denote $p$ repeated
opening parentheses by $(^p$ and $p$ repeated opening parentheses
followed by $p$ closing parentheses as $(\ldots)_p =
(^p\ldots)^p$. Furthermore, we use the special notation $[s,t,p]$ for
matchings of the form $(()_s()_t)_p$.
\begin{conj}
\label{co:asm2hex}
The coefficient $a_F$ of matching $F=[s,t,p]$ in the stationary state
of $H^{\rm C}_{2p+2s+2t}$ is given by
\begin{eqnarray*}
a_{[s,t,p]} &=& \det_{1\leq i,j\leq s} \left( \binom{2(s+t+p)+j-2i}{s+t-j} -
\binom{2(s+t+p)+j-2i}{s+t-j-2i+1} \right) \\
&=& \prod_{j=1}^s \frac{(j-1)!(2t+2p+2j-1)!(2p+2j)_j(3t+2p+3j)_{s-j}}
{(t+2p+s+2j-1)!(t+s-j)!}
\end{eqnarray*}
\end{conj}
The evaluation of the determinant in this conjecture is due to Ciucu
and Krattenthaler \cite{CiucuK02}.

The numbers appearing in Conjecture \ref{co:asm2hex} enumerate lozenge
tilings, or dimer configurations, on hexagons $H_{\rm
a}(s,2p+s+t-1,t)$ with maximal staircases removed from adjacent
vertices (see \cite{CiucuK02} for definitions). With Conjecture
\ref{co:TL-FPL} in mind, Conjecture \ref{co:asm2hex} implies 
that FPL diagrams with matching $[s,t,p]$ are equinumerous with
lozenge tilings on $H_{\rm a}(s,2p+s+t-1,t)$. Here we give the
explicit bijection.
%
\begin{theorem}
\label{th:Ha}
Let $n=p+s+t$. Then
\[
M_{[s,t,p]}(G^{\rm V}_{2n}) = \prod_{j=1}^s
\frac{(j-1)!(2t+2p+2j-1)!(2p+2j)_j(3t+2p+3j)_{s-j}}
{(t+2p+s+2j-1)!(t+s-j)!}
\]
\end{theorem}

In the rest of this section we give a sketch of the proof of Theorem
\ref{th:Ha}. We first note that every matching fixes the loops of the
corresponding FPL diagrams to pass through a particular set of
edges. Let us call these edges {\bf fixed edges}. The way these edges
are fixed is prescribed in 
\begin{lemma}
The implication
\[
\begin{picture}(140,30)
\put(0,0){\epsfxsize=45pt\epsfbox{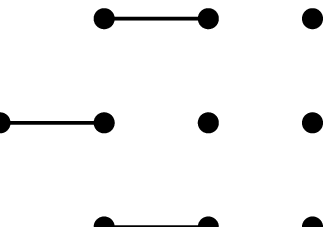}}
\put(65,10){$\Rightarrow$}
\put(95,0){\epsfxsize=45pt\epsfbox{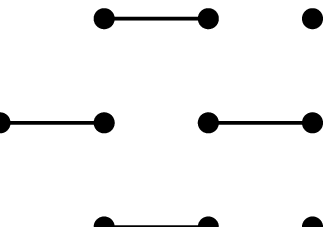}}
\end{picture}
\]
holds if the top and bottom loop segments do not belong to the same
loop and if either of the following holds,
\begin{itemize}
\item[i)]
The middle loop segment belongs to a third loop.
\item[ii)] 
If the middle loop segment belongs to the same loop as the top or
bottom segment, it is connected to it via one of the leftmost
edges. 
\end{itemize}
\end{lemma}
In case i) the third loop has to pass between the top and bottom loop
segments while the condition in case ii) excludes possibilities like 
\[
\begin{picture}(45,30)
\put(0,0){\epsfxsize=45pt\epsfbox{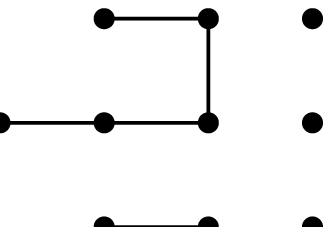}}
\end{picture}
\]

\begin{remark}
The fixed edges can be adjacent and therefore do not necessarily form
a perfect matching. Furthermore, for a matching $[s,t,p]$ every site
of $G^{\rm V}_{2n}$ is part of a fixed edge. 
\end{remark}

\begin{remark}
None of the FPL diagrams with matching $[s,t,p]$ contains an internal
closed loop.
\end{remark}

\begin{example}
For the matching $(2,1,3)=((((())())))$ the fixed edges resulting from
the first $5$ parentheses, and the parentheses $6,7$ and $8$ are
\[
\begin{picture}(200,110)
\put(0,0){\epsfxsize=200pt\epsfbox{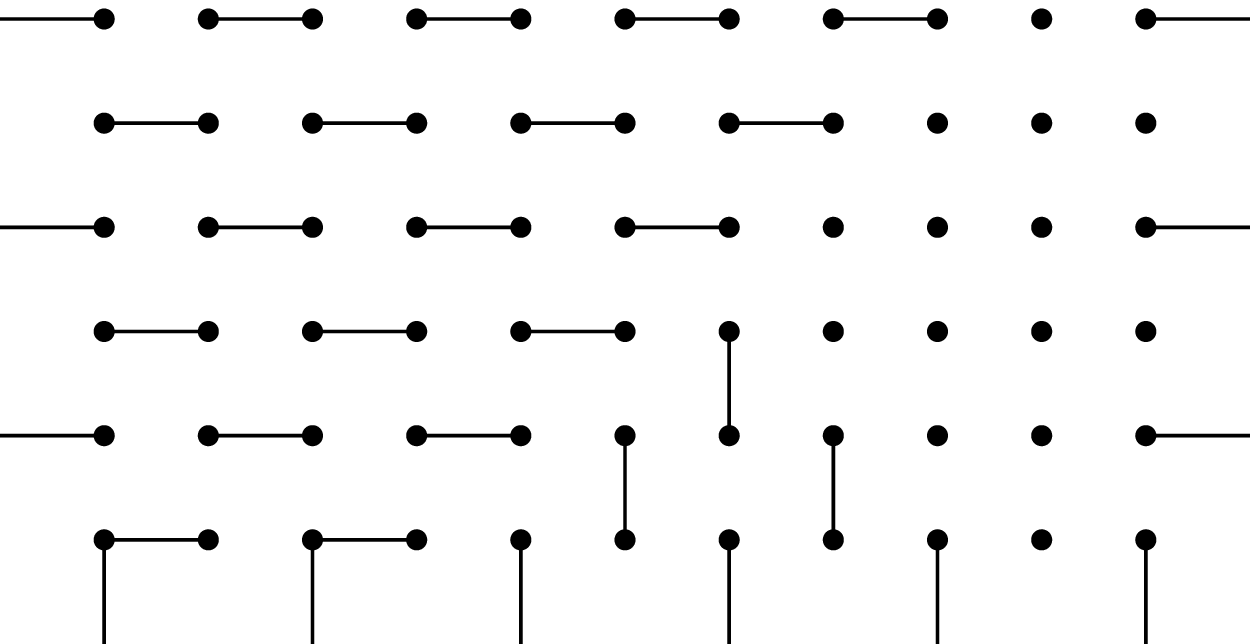}}
\end{picture}
\]
For clarity of the picture, the fixed edges resulting from the last $4$
parentheses are not drawn. That region of fixed edges partly overlaps
the region of fixed edges corresponding to the first $5$ parentheses.
\end{example}

\begin{remark}
\label{rem:fixed}
More edges may in fact be fixed due to the specific geometry of the grid
$G^{\rm V}_{2n}$. 
\end{remark}

Since for matchings $[s,t,p]$ all sites correspond to a fixed edge,
the problem of counting all FPL diagrams with this matching reduces to
a dimer problem which is equivalent to a lozenge tiling problem on a
graph defined by the fixed edges. This graph is called the {\bf fixed
graph}. It forms a folded piece of triangular lattice and is
obtained by joining the midpoints of all the fixed edges. 

\begin{example}
\label{ex:fixedG}
For the matching $(2,1,3)$, the set of fixed edges and
its corresponding fixed graph is
\[
\begin{picture}(200,135)
\put(0,0){\epsfxsize=200pt\epsfbox{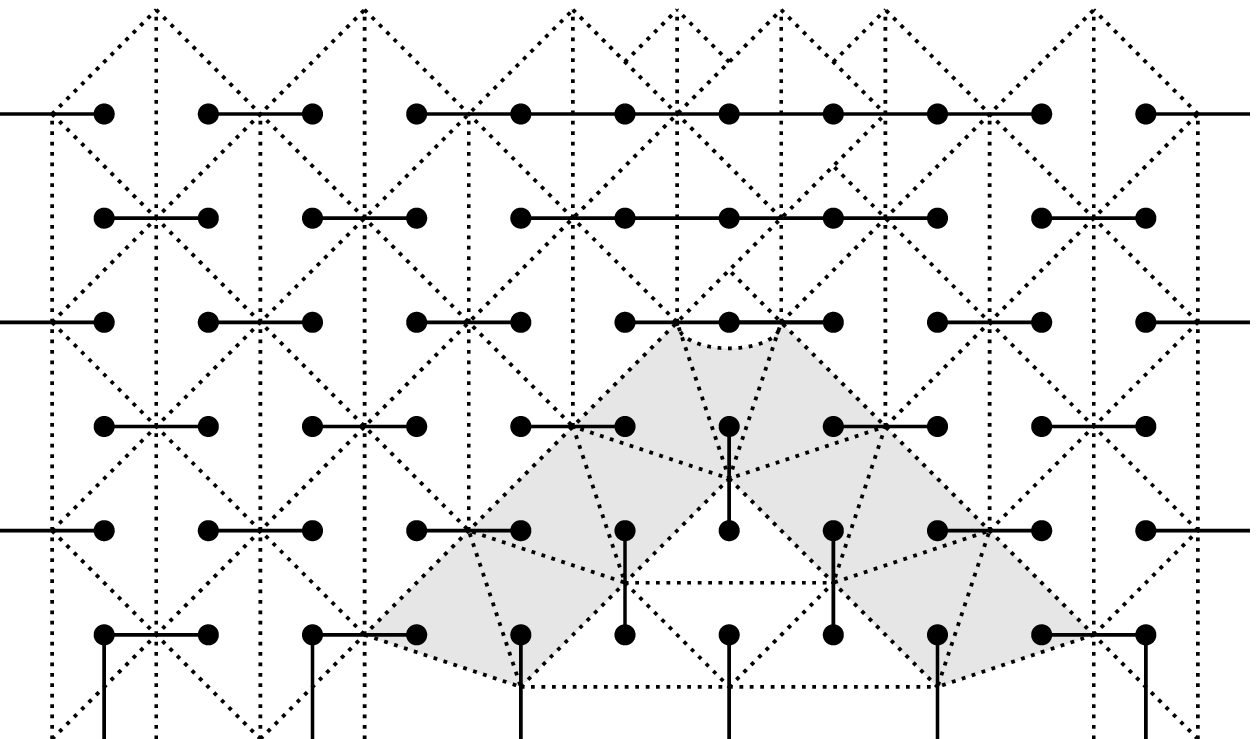}}
\end{picture}
\]
The graph consists of three homogeneous regions arising from the fixed
edges of the first 5 parentheses, parentheses 6,7 and 8 and the last 4
parentheses respectively. The first and last region partially overlap
in the top middle of the picture, where the fixed edges are
adjacent. To guide the eye, the region between the three homogeneous
regions is shaded. 
\end{example}

As already hinted at in Remark \ref{rem:fixed}, the particular geometry of
the fixed graph completely determines dimer or tiling configurations
on certain regions. These regions are the overlapping parts of the
fixed graph and patches in the lower left and right corners. They
do therefore not contribute to the enumeration and may as well
be removed. 

\begin{example}
By removing the overlapping part of the graph in Example
\ref{ex:fixedG} and deforming it so that it fits on the 
triangular lattice, it follows that the fixed graph for the matching
$(2,1,3)$ is 
\[
\begin{picture}(200,90)
\put(0,0){\epsfxsize=200pt\epsfbox{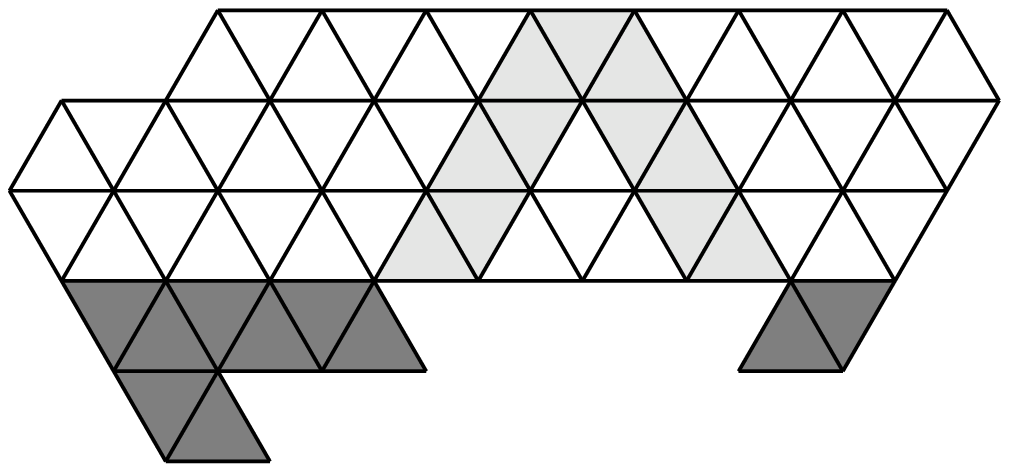}}
\end{picture}
\]
The lightly shaded area corresponds to the shaded area of Example
\ref{ex:fixedG}, and is a guide to the eye only. The regions on which
dimer configurations are enforced by the geometry are shaded darker. 
\end{example}

When the regions with fixed dimer configurations are removed, the
resulting graph is precisely that of a hexagon with maximal staircases
removed from adjacent corners \cite{CiucuK02} and whose sides are given
by $s$, $2p+s+t-1$ and $t$. The number of dimer configurations on such
hexagons is given by Theorem 1.6 of \cite{CiucuK02}. Theorem
\ref{th:Ha} then follows immediately from that result.

\section{Concluding remarks}
There is another model known that displays the appearance of
alternating sign matrix numbers in its stationary state. This is the
rotor model \cite{BatchGN02} which is based on two Temperley-Lieb
algebras. It is conjectured that in that model also 3-enumerations of ASMs
\cite{Kupe00} play a role. 

Many conjectures are stated but remain unproven. One may
hope that again the solvability of the six-vertex model and the XXZ spin
chain \cite{Baxter82} can be used to find proofs. In particular
the very special properties at $q=\exp (\i \pi/3)$, see
e.g. \cite{Baxter89}, have already led to several interesting 
results \cite{FridSZ00,FridSZ01,GierBNM02,KitaMST02,Strog01,Strog02}.

\section*{Acknowledgment}
This work has been supported by the Australian Research 
Council (ARC). I thank Murray Batchelor, Saibal Mitra, Bernard
Nienhuis, Paul Pearce and Vladimir Rittenberg for stimulating and
useful discussions and for collaborations where parts of the results
presented here were obtained. I am furthermore grateful to Christian
Krattenthaler and an anonymous referee for very useful suggestions, and
lastly thank Catherine Greenhill and Ole Warnaar for a critical
reading of the manuscript.

\end{document}